\documentclass[a4paper,12pt,reqno]{amsart}
\usepackage{vmargin}
\usepackage[pdftex]{graphicx}
\usepackage{amsfonts}
\usepackage{amssymb}
\usepackage{mathrsfs}
\usepackage{stmaryrd}
\usepackage{amsmath}
\usepackage{amsthm}
\usepackage[shortlabels]{enumitem}
\usepackage{nomencl}
\usepackage[bookmarks=true]{hyperref}   
\usepackage{esint}
\usepackage{dsfont}
\usepackage{upgreek}
\usepackage{xcolor}
\usepackage{slashed}
\usepackage{scalerel}
\usepackage{comment}
\usepackage{todonotes}
\usepackage[utf8]{inputenc}

\makeatletter
\renewcommand\section{\@startsection{section}{1}{\z@}%
                       {-3\p@ \@plus -4\p@ \@minus -4\p@}%
                       {3\p@ \@plus 4\p@ \@minus 4\p@}%
                      {\normalfont\normalsize\centering\scshape}}
\makeatother

\numberwithin{equation}{section} 	

\hypersetup{
    colorlinks,%
}

\author{Lashi Bandara}
\title[Fun. calc. and harm. anal. in geometry]
{Functional calculus and harmonic analysis in geometry}

\date{\today}

\address{Lashi Bandara, 
Institut für Mathematik,
Universität Potsdam, 
D-14476, Potsdam OT Golm, Germany
}
\urladdr{\href{http://www.math.uni-potsdam.de/~bandara}{http://www.math.uni-potsdam.de/~bandara}}
\email{\href{mailto:lashi.bandara@uni-potsdam.de}{lashi.bandara@uni-potsdam.de}}



\def\colour{\colour}

\def\colour{\color}

\newtheorem{theorem}{Theorem}[section]

\newtheorem{lemma}[theorem]{Lemma}

\newcommand{\cbrac}[1]{\left(#1\right)}

\newcommand{\dbrac}[1]{\left\{#1\right\}}
\newcommand{\modulus}[1]{|#1|}

\newcommand{\set}[1]{\dbrac{#1}}

\newcommand{\dom}{ {\mathrm{dom}}}
\newcommand{\ran}{ {\mathrm{ran}}}
\newcommand{\nul}{ {\mathrm{ker}}}

\newcommand{\comp}{\, \circ\, }

\newcommand{\e}{\mathrm{e}}

\newcommand{\R}{\mathbb{R}}
\newcommand{\C}{\mathbb{C}}
 
\newcommand{\In}{\mathbb{Z}}

\newcommand{\Na}{\ensuremath{\mathbb{N}}}

\newcommand{\script}[1]{\mathscr{#1}}

\renewcommand{\emptyset}{\varnothing}
\newcommand{\union}{\cup}

\newcommand{\intersect}{\cap}

\newcommand{\disunion}{\sqcup}

\newcommand{\rest}[1]{{{\lvert_{}}_{}}_{#1}}
\newcommand{\close}[1]{\overline{#1}}		

\newcommand{\ind}[1]{\raisebox{\depth}{\(\chi\)}_{#1}}	

\newcommand{\id}{{\rm id}}
\renewcommand{\epsilon}{\varepsilon}

\renewcommand{\phi}{\varphi}

\newcommand{\ch}[1]{\upchi_{{#1}}}	

\newcommand{\embed}{\hookrightarrow}		

\newcommand{\tensor}{\otimes}

\newcommand{\norm}[1]{\| #1 \|}			
\newcommand{\snorm}[1]{\left\| #1 \right\|}			

\newcommand{\spt}[1]{{\rm spt} {\text{ }}#1}	

\DeclareMathOperator{\gap}{\hat{\updelta}}		

\DeclareMathOperator{\diam}{diam}		

\DeclareMathOperator{\len}{\ell}			

\DeclareMathOperator{\divv}{div}		

\newcommand{\Ric}{{\rm Ric}}			
\newcommand{\SecondFF}{\mathrm{II}}			

\DeclareMathOperator{\inj}{inj} 		

\newcommand{\bnd}{\partial}			
\newcommand{\interior}[1]{\mathring{#1}}	
\DeclareMathOperator{\proj}{\mathbf{P}}	
\newcommand{\tanb}{{\rm T}}		
\newcommand{\cotanb}{{\rm T}^\ast}	

\newcommand{\on}{\vec{\mathrm{n}}}

\DeclareFontFamily{OT1}{restrictfont}{}
\DeclareFontShape{OT1}{restrictfont}{m}{n}{<-> fmvr8x}{}

\newcommand{\adj}[1]{{#1}^\ast}			
\newcommand{\extd}{{\rm d}}			

\newcommand{\inprod}[1]{\left\langle #1 \right\rangle}	
\newcommand{\conn}[1][{}]{{\nabla_{{#1}}}}		
\newcommand{\Leb}[1][{}]{\script{L}^{#1}}			

\DeclareMathOperator{\Spinors}{\slashed{\Delta}}	
\newcommand{\spin}[1]{\slashed{#1}}		

\newcommand{\sym}{\upsigma}

\newcommand{\spec}{\mathrm{spec}}		

\newcommand{\Lp}[2][{}]{{\rm L}^{#2}_{\rm #1}}		
\newcommand{\Ck}[2][{}]{{\rm C}^{#2}_{\rm #1}}		
\newcommand{\Hard}[2][{}]{{\rm H}^{#2}_{\rm #1}}		
\newcommand{\SobH}[2][{}]{\Hard[#1]{\rm #2}}
\newcommand{\Lips}[1][{}]{{\rm Lip}_{\rm #1}}		

\newcommand{\Sec}[1]{\mathrm{S}_{#1}}
\newcommand{\OSec}[1]{\mathrm{S}^\mathrm{o}_{#1}}

\newcommand{\maxx}[1]{\left<#1\right>}

\renewcommand{\Re}{\mathrm{Re}\ }

\newcommand{\iden}{{\mathrm{I}}}
\DeclareMathOperator{\sgn}{sgn}

\newcommand{\Hil}{\script{H}}			
\newcommand{\Lap}{\Delta}			

\newcommand{\Q}[1][{}]{Q_{#1}}			
\newcommand{\DyQ}{\script{Q}}			

\newcommand{\scale}{\mathrm{t_S}}			

\newcommand{\sE}{\script{E}}

\newcommand{\cB}{\mathcal{B}}

\newcommand{\cE}{\mathcal{E}}
\newcommand{\cF}{\mathcal{F}}

\newcommand{\cM}{\mathcal{M}} 
\newcommand{\cN}{\mathcal{N}}

\newcommand{\cW}{\mathcal{W}}

\newcommand{\mg}{\mathrm{g}}
\newcommand{\mgt}{\tilde{\mg}}
\newcommand{\mh}{\mathrm{h}}

\newcommand{\Av}{\mathbb{E}}
\newcommand{\Pri}{\upgamma}

\newcommand{\Div}{\mathrm{L}}

\newcommand{\Carl}{\mathcal{C}}

\newcommand{\CBox}{\mathrm{R}}

\newcommand{\U}{\mathrm{U}}

\newcommand{\Sph}{\mathrm{S}}

\newcommand{\Dir}{{\rm D} }
\newcommand{\DirA}{{\rm A} }

\newcommand{\BDir}{\slashed{\DirA}}

\newcommand{\cc}{\mathrm{cc}}

\newcommand{\Ppb}{{\rm P}}

\newcommand{\QQ}{{\rm\bf Q}}
\newcommand{\PP}{{\rm\bf P}}

\newcommand{\dtt}[1][{t}]{\frac{d#1}{#1}}

\newcommand{\met}{\uprho}		

\newcommand{\Ball}{\mathrm{B}}

\newcommand{\checkH}{\check{\mathrm{H}}}	

\newcommand{\SDir}{\slashed{\Dir}}
\newcommand{\SDirRS}{\Dir_{\mathrm{RS}}}

\newcommand{\hk}{\uprho}

\newcommand{\Hinfty}{\mathrm{H}^\infty}

\newcommand{\RNum}[1]{\uppercase\expandafter{\romannumeral #1\relax}}

\begin{document}

\maketitle
\vspace*{-2em}
\begin{abstract}
In this short survey article, we showcase a number of non-trivial geometric problems that have recently been resolved by marrying methods from functional calculus and real-variable harmonic analysis.
We give a brief description of these methods as well as their interplay.  
This  is a succinct survey that hopes to  inspire geometers and analysts alike to study these methods so that they can be further developed to be potentially applied to a broader range of questions.
\end{abstract}
\tableofcontents
\vspace*{-2em}

\parindent0cm
\setlength{\parskip}{\baselineskip}
\setlist{noitemsep, topsep=-.5\parskip,leftmargin=*, listparindent=0pt}

\renewcommand{\nomname}{Notation}
\makenomenclature

\section{Introduction}

The purpose of this short survey article is to highlight some geometric questions that have recently been addressed via the use of functional calculus and harmonic analytic methods. 
In the literature, these phrases are umbrella terms that carry different meanings in different schools. 
Therefore, before we embark on considering the details, let us first note what we mean by these terms.

Functional calculus  is the ability to take functions of operators and manipulate them as if they were functions.
This is an extraordinarily powerful way to think and one motivation is that it provides a conceptual method to solve partial differential equations, particularly evolution equations. 
Harmonic analysis, on the other hand, is the art in which a mathematical object can be seen as a ``signal'', and whose goal is to decompose this signal, typically in some scale invariant way, to simpler parts which are mathematically tractable.
Perhaps the simplest of this kind is the Fourier series, where every square integrable function on the circle can be written as an infinite sum of trigonometric functions, where the properties of these latter functions are well known and understood.
Geometry, to us, almost always means the shape of space.
In practice, this is quantified in a multitude of ways, but here, we simply contrast it from the topology or structure of the space, which is seldom of concern to us.

At a first glance, perhaps it might seem peculiar to talk about a combination of functional calculus and harmonic analysis.
However, this is an extremely fruitful marriage, and  the late Alan McIntosh  was one of the pioneers of this school of thought. 
His idea  was to demonstrate the ability to construct a powerful and desirable functional calculus, namely the \emph{$\Hinfty$ functional calculus}, via \emph{quadratic estimates}. 
These estimates can be interpreted as the ability to \emph{reconstruct a signal in norm}  via operator adapted \emph{band-pass filters}, where these filters only rely on the spectral properties of the operator.
Taking a perspective of Fourier theory as a functional calculus for the Laplacian, the $\Hinfty$ functional calculus can be seen as a Fourier theory adapted to an operator.
The quadratic estimates point of view of this functional calculus is extremely useful  as it provides a bridge to  methods emerging from real-variable harmonic analysis in order to compute these estimates.

An important motivating factor for the development of these methods was the treatment of the \emph{Kato square root problem}, which remained open for almost half a century.
In \cite{KatoSq}, Kato formulated the most abstract version of the problem. 
He first demonstrated that if $A$ is a maximally-accretive operator, meaning that the spectrum and numerical range sit in the complex plane with positive real part,  then for any $\alpha < \frac{1}{2}$, the domains $\dom(A^\alpha) = \dom({\adj{A}}^\alpha)$.
The equality was known to be invalid in general when $\alpha > \frac{1}{2}$ and Lions in \cite{Lions} showed that that this equality even fails in general for maximally-accretive operators for the critical case $\alpha = \frac{1}{2}$.
However, the question remained open as to whether the critical case was true for regularly-accretive operators, which are operators whose numerical range and spectrum are contained in a proper sector of the complex plane with positive real part.
In 1972, McIntosh authored \cite{Mc72} which  dealt a fatal blow to the abstract question by constructing a regularly-accretive operator $A$ for which $\dom(\sqrt{A}) \neq \dom(\sqrt{\adj{A}})$.
Despite this fact, given that the main interest of Kato's original question was  in its applications to partial differential equations consisting of operators in divergence form, the Kato square root problem as we know it today became rephrased as to determine  whether
$$ \dom(\sqrt{-\divv B \nabla}) = \SobH{1}(\R^n) = \set{u \in \Lp{2}: \nabla u \in \Lp{2}}$$
for coefficients $B \in \Lp{\infty}(\R^n)$ with a constant $\kappa > 0$ such that 
$$\Re \inprod{B(x) u, u}_{\R^n} \geq \kappa \modulus{u}^2$$ 
almost-everywhere in $\R^n$. 
If the coefficients $B(x)$ are real and symmetric, then the desired equality follows simply from operator theoretic considerations alone.
It is the complex, non-symmetric coefficient case that is of significance. 

The resolution of this question for $n = 1$ was the famous paper \cite{CoMcM} due to Coifman,  McIntosh and Meyer which they authored in  1981.
It is a paper that has spawned a number of new directions in mathematics including wavelet theory. 
Moreover, this paper made significant contributions to the theory of singular integrals as well as the development of the so-called $T(b)$ theorems.
Beyond their general usefulness, this was of paramount importance to the resolution of the Kato square root problem in arbitrary dimensions.

In one dimension, the Kato square root problem is actually equivalent to the boundedness of the Cauchy integral operator on a Lipschitz curve.
However, it is the description of this problem from a differential operator point of view, where the functional calculus perspective takes centre stage, that was a key observation to its resolution. 
The general case was finally settled in 2002 by Auscher, Hofmann, Lacey, McIntosh and Tchamitchian in \cite{AHLMcT}.
There are numerous surveys of the Kato square root problem and its resolution.
Namely, the lecture series \cite{Auscher} by Auscher, as well as the survey article \cite{HMc} by Hofmann and McIntosh give detailed accounts.
See also references therein.

The general Kato square root problem was initially proved using second-order methods in \cite{AHLMcT}.
In 2005, in \cite{AKMc}, Axelsson (Rosén), Keith and McIntosh produced a first-order perspective of this problem. 
There, they constructed  a larger first-order Dirac-type operator which encoded the problem. 
The idea  was to rewrite  the divergence form equation into the system
$$ \Pi_{B} = \begin{pmatrix} 0 & -\divv B \\
			 \nabla & 0
		\end{pmatrix}$$
and by showing that this operator has an $\Hinfty$ functional calculus, the desired conclusion was obtained. 
The required domain equality follows from showing that the functions 
$$ \zeta \mapsto \frac{\sqrt{\zeta^2}}{\zeta}\quad \text{and} \quad \zeta \mapsto \frac{\zeta}{\sqrt{\zeta^2}}$$
applied to $\Pi_{B}$ define a bounded operator.
In \cite{AKMc}, the authors also demonstrate how to prove a similar problem for perturbations of the Hodge-Dirac operator of the form $D_B = \extd + B^{-1} \adj{\extd} B$. 
They further obtain corresponding results on compact manifolds.
In fact, it  this paper which provided the key conceptual point of view to make similar methods work in more geometric settings. 
It is important to also acknowledge Morris, where in his thesis \cite{MorrisThesis}, he developed many necessary technical results on metric spaces which were later adapted  to attack geometric problems on vector bundles over manifolds.

Despite the fact that these methods were developed primarily for non-smooth Euclidean problems, their success in analysing differential operators on vector bundles has been for reasons similar in spirit to their success in the non-smooth setting. 
Pseudo-differential or Fourier analytic methods often rely on localisation and smoothness.
The functional calculus that we will talk about here is sufficiently implicit and general to be able to treat non-smooth coefficients, which requires a global point of view. 
In a similar vein, bundles typically carry non-trivial geometry, and even if the coefficients of the problem are smooth, the tools need to be sufficiently abstract in order to capture the global nature without resorting to local descriptions.
This is really the key indication that the methods borne out of the resolution of the Kato square root problem, and in particular its first-order reincarnation,  might have useful applications to geometry, in particular to analysis on bundles.

Lastly, we  remark that this survey is far from being comprehensive.
The results we discuss here are almost exclusively contributions of the author and co-authors.
It represents only a small subset of the successes of coupling functional calculus and harmonic analysis to be an effective force in solving problems in geometry. 
The author does not have sufficient expertise to comment on the general picture, and so the limited focus here is for the sake of conciseness and accuracy. 
Moreover, the author apologises in advance for the many omitted names and references.
It is the author's hope that geometers and analysts become better aware of these methods and  embark on further developing these tools.
In particular, there is a real need to adapt the methods arising from real-variable harmonic analysis, which are typically Euclidean in nature, to better account for the underlying geometry.
\section*{Acknowledgements}
The author was supported by SPP2026 from the German Research Foundation (DFG).
David Rule deserves mention as this paper was borne out of a Colloquium talk given at Linköping University upon his invitation.
Andreas Rosén needs to acknowledged for useful feedback on an earlier version of this article. 
The anonymous referee also deserves a mention for their helpful suggestions.
Furthermore, the author acknowledges the gracious hospitality of his auntie Harshya Perera in Sri Lanka who hosted him during the time in which this paper was written. 
\section{Functional Calculus}
\label{Sec:FunCalc} 

\subsection{Fourier series}
Let us begin our exposition by stating a beautiful and well known fact seen from a seldom mentioned perspective.

Recall that for a function $u \in \Lp{2}(\Sph^1)$, we can write its Fourier series as 
$$u(\theta) = \sum_{n=-\infty}^\infty a_n \e^{\imath n \theta},$$
where the coefficients $a_n \in \C$ are uniquely determined by $u$.
Now, consider the Laplacian $\Lap_{\Sph^1}$ on $\Lp{2}(\Sph^1)$, or equivalently, $-\frac{d^2}{d\theta^2}$ on $[0,1]$ with periodic boundary conditions on $\Lp{2}([0,1])$.
It is well known that the spectrum of this operator consists only of eigenvalues and it is precisely $\spec(\Lap_{\Sph^1}) = \set{\lambda_n = n^2}_{n=0}^\infty$. 
The corresponding eigenfunctions are $\set{\e^{\imath n \theta}, \e^{- \imath n \theta}}_{n=0}^\infty$.
\nomenclature{$\spec(T)$}{Spectrum of an operator $T$.}
When $u \in \dom(\Lap_{\Sph^1})$, the operator $\Lap_{\Sph^1}$ acts on $u$ as
\begin{equation}
\label{Eq:SphAct} 
(\Lap_{\Sph^1}u) (\theta) = \sum_{n=-\infty}^\infty n^2 a_n \e^{\imath n \theta}.
\end{equation}
Therefore, we see the decomposition of $u$ as a Fourier series is, in fact, the decomposition of $f$ with respect to the frequencies of the Laplacian on $\Sph^1$. 

Via \eqref{Eq:SphAct}, we can construct a rudimentary \emph{functional calculus} for $\Lap_{\Sph^1}$. 
We call a function $f:\In \to \C$ a \emph{symbol} and define a function of $\Lap_{\Sph^1}$ by $f$ via: 
$$ (f(\Lap_{\Sph^1})u) (\theta) := \sum_{n=-\infty}^\infty f(n^2)a_n \e^{\imath n \theta}.$$
Here, the significance of spectral theory  in the construction of this functional calculus is explicit. 
Indeed, this will be true for functional calculi we consider for more general operators.
For an arbitrary function $f$, we would need to compute the domain of this operator, which is precisely the $u \in \Lp{2}(\Sph^1)$ for which the  defining series converges. 
Moreover, it is easy to see that if $f$ is bounded, then $f(\Lap_{\Sph^1})$ is a bounded operator on $\Lp{2}(\Sph^1)$.

Note that this functional calculus, at least for bounded functions, enjoys a few desirable properties. Of significance is that it is an algebra homomorphism $f \mapsto f(\Lap_{\Sph^1})$, which satisfies 
$$fg \mapsto f(\Lap_{\Sph^1})g(\Lap_{\Sph^1}), \quad (f\comp g) \mapsto f( g(\Lap_{\Sph^1})),\quad\text{and}\quad 1 \mapsto \iden.$$
These are essential features of the functional calculi we consider in this survey. 

An important motivation for  functional calculi comes from evolution equations, where solutions can be simply generated and represented through the functional calculus.
For example, let us consider the following heat equation:
\begin{equation*}
\partial_t u(t, \theta) +  \Lap_{\Sph^1} u(t, \theta) = 0, \qquad \lim_{t \to 0} u(t, \theta) = u_0 (\theta). 
\end{equation*}
The unique solution to this is given by $u(t,\theta) =  (\e^{-t \Lap_{\Sph^1}}u_0)(\theta) = \sum_{n=-\infty}^\infty e^{-t n^2} a_n e^{\imath n \theta}$.

\subsection{Non-negative self-adjoint operators with discrete spectrum}
\label{Sec:Discrete}
Taking inspiration from this very concrete and classical example, let us instead consider the following abstract setup. 
Let $\Hil$ be a Hilbert space, $\Dir$ a non-negative self-adjoint operator on $\Hil$ with discrete spectrum $\spec(\Dir) = \set{\lambda_i \geq 0}$.
We denote its eigenfunctions by $\set{\psi_i}$, which are orthogonal due to self-adjointness and as before, whenever $u \in \Hil$, we obtain coefficients $a_n \in \C$ such that $u = \sum_{n} a_n \psi_n$. 
In this context, a functional calculus for $\Dir$ can be generalised \emph{mutatis mutandis} from the setting of $\Sph^1$ in the following way: 
\begin{equation*} 
f(\Dir) u := \sum_{n} f(\lambda_n) a_n \psi_n.
\end{equation*}
As before, the $f$ needs to be chosen from an appropriate class in order to obtain a well-defined expression and to understand the mapping properties of $f(\Dir)$.
For instance, it is easy to see that if $f: \R \to \C$ is bounded and continuous, the operator $f(\Dir)$ is well-defined and bounded on $\Hil$. 

There are good reasons for considering this more abstract setup as the following example illustrates.
Let $(\cM,\mg)$ be a smooth compact Riemannian manifold (without boundary), $\Hil = \Lp{2}(\cM) = \Lp{2}(\cM,\mg)$, and $\Dir = \Lap_{\cM,\mg} = \nabla^{\ast,g} \nabla$, the Laplacian on $(\cM,\mg)$. 
Then, $\id: \SobH{1}(\cM) \embed \Lp{2}(\cM)$ is a compact map and $\dom(\Lap_{\mg}) = \SobH{2}(\cM) \subset \SobH{1}(\cM)$.
The combination of these two facts yield that $\spec(\Lap_{\cM})$ is discrete with no finite accumulation points. 
This exactly fits the abstract picture we have just painted.
Of particular significance is the fact that we have not alluded to the exact value of the eigenvalues and  their corresponding eigenfunctions, since unlike $\Lap_{\Sph^1}$ on $\Sph^1$, in general,  they cannot  be computed.

It is a maxim that the better you know the spectrum, the larger the class of functions that you can consider in order to construct a functional calculus. 
Since $\Dir$ is a self-adjoint operator, it further admits a functional calculus for \emph{Borel} functions $f: \R \to \C$, which is obtained by extracting a spectral measure. We shall not indulge in the details of its construction but it can be found in \cite{Kato} by Kato. 

We now touch upon an important property of these functional calculi - the  \emph{reconstruction of a signal in norm}.
Let us consider $u \in \Hil$ as a ``signal''.
We want to use the functional calculus of $\Dir$ to reconstruct $\norm{u}$ but up to a constant.
To keep our exposition simple, let us assume $\nul(\Dir) = 0$. 
This is not a severe restriction: if $\nul(\Dir) \neq 0$, then we can split $\Hil = \nul(\Dir) \oplus^\perp \close{\ran(\Dir)}$ and instead work on $\close{\ran(\Dir)}$ as the new space.
Let $\phi: \R \to \R$ be a piecewise smooth function such that $\psi \neq 0$ and for which there exists $\alpha > 0$ and $C > 0$ satisfying:
\begin{equation}
\label{Eq:Bpassdec} 
\modulus{\psi(x)} \leq C \min \set{\modulus{x}^\alpha, \modulus{x}^{-\alpha}}
\end{equation}
for almost-every $x \in \R$.
We want to justify the idea that $\psi(t\Dir)$ is a \emph{band-pass filter} for frequencies localised around $\frac{1}{t}$.
Let us first compute with $\psi(x) = \chi_{[1,2]}(x)$ in order to make the exposition clear.
\begin{figure}[ht]
\begin{centering}
\includegraphics{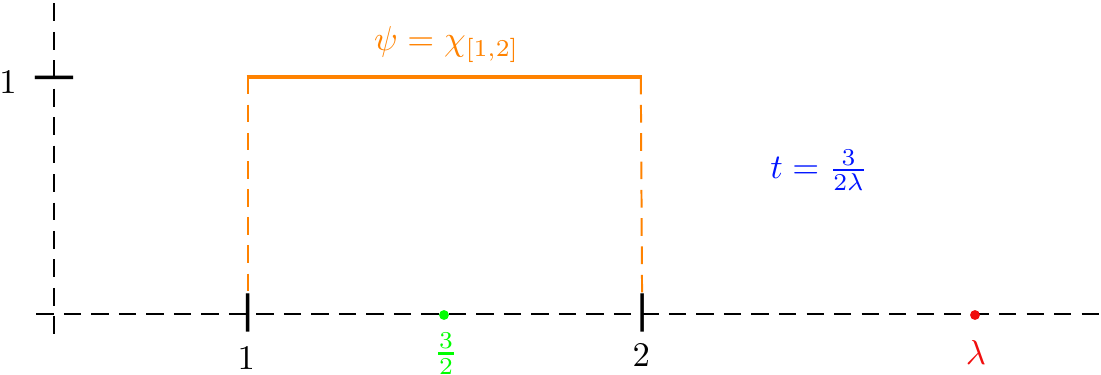}
\end{centering}
\end{figure}

On writing $u = \sum_{n} a_n \psi_n$, note that 
\begin{align*}
\norm{\chi_{[1,2]}(t\Dir)u}^2 
	&= \sum_{n} a_n^2 \norm{\chi_{[1,2]}(t\Dir)\psi_n}^2   \\ 
	&= \sum_n a_n^2 \norm{\chi_{[1,2]}(t\lambda_n) \psi_n}^2 
	= \sum_n a_n^2 \modulus{\chi_{[1,2]}(t\lambda_n)}^2 \norm{\psi_n}^2.
\end{align*}
Therefore,
\begin{equation}
\label{Eq:2}
\int_0^\infty \norm{\chi_{[1,2]}(t\Dir)\psi_n}^2\ \frac{dt}{t} = \norm{\psi_n}^2 \int_0^\infty \modulus{\chi_{[1,2]}(t\lambda_n)}^2\ \frac{dt}{t}.
\end{equation}
On substituting $s = t \lambda_n$ so that $dt = \frac{ds}{\lambda_n}$ and $\frac{1}{t} = \frac{\lambda_n}{s}$, we obtain 
\begin{equation}
\label{Eq:Haar} 
\begin{aligned}
\int_0^\infty \modulus{\chi_{[1,2]}(t \lambda_n)}^2\ \frac{dt}{t} 
	&= \int_{0}^\infty \modulus{\chi_{[1,2]}(s)} \frac{ds}{\lambda_n } \cdot \frac{\lambda_n}{s}	  \\  
	&= \int_0^\infty \modulus{\chi_{[1,2](s)}}^2\ \frac{ds}{s}  
	= \int_1^2 \frac{1}{s}\ ds  = \log(2). 
\end{aligned}
\end{equation} 
Combining \eqref{Eq:2} with \eqref{Eq:Haar}, we obtain
$$\int_0^\infty \norm{\chi_{[1,2]}(t\Dir)u}^2\ \frac{dt}{t} = \sum_{n} a_n \log(2) \norm{u_n}^2 = \log(2) \norm{u}^2.$$

If we were to repeat this exercise with a general $\psi$ satisfying \eqref{Eq:Bpassdec}, then by a similar calculation, we would obtain that
\begin{equation*}
\int_0^\infty \norm{\psi(t\Dir)u_n}^2\ \frac{dt}{t} = \norm{u_n}^2 \int_0^\infty \modulus{\psi(t\lambda_n)}^2\ \frac{dt}{t} \simeq \norm{u_n}^2
\end{equation*}
so that 
\begin{equation*}
\int_0^\infty \norm{\psi(t\Dir)u}^2\ \frac{dt}{t} \simeq \norm{u}^2.
\end{equation*}

We emphasise here that in \eqref{Eq:Haar}, the cancellation of the spectral points $\lambda_n$ under a change of variables occurs precisely due to the fact that we compute with respect to the Haar measure $\frac{dt}{t}$ on $\R$.
It is instructive to repeat this simple calculation with respect to another measure to better understand the significance of $\frac{dt}{t}$. 

\subsection{Bisectorial operators and the McIntosh Theorem}
The key idea in the calculation in the previous section is that  $\psi(t \Dir)$ filters out all frequencies but those concentrated near $\frac{1}{t}$, and at least for the class of invertible self-adjoint operators with discrete spectrum that we considered, and despite the lack of ability to explicitly compute the spectrum,  we were able to reconstruct the signal in norm up to a constant. 
We will see here that, in fact, this is a characterising condition for the existence of a more significant and useful functional calculus that can be applied to a broader class of operators. 

Before we introduce this class of operators, for $\alpha < \pi/2$,  let us define 
$$\OSec{\alpha} := \set{\zeta \in \C\setminus\set{0}: \modulus{\arg(\zeta)} <  \alpha\quad\text{or}\quad \modulus{\arg(-\zeta)} < \alpha},$$ 
the open bisector of angle $\alpha$ in the complex plane, and $\Sec{\alpha} := \close{\OSec{\alpha}}$, the closed bisector of angle $\alpha$.
For $\omega < \pi/2$, an operator $T: \dom(T) \subset \Hil \to \Hil$ is called $\omega$-bisectorial if: 
\begin{enumerate}[(i)]
\item $T$ is closed, 
\item $\spec(T) \subset \Sec{\omega}$, and
\item 
\label{Def:Sec:RE} 
for all $\mu \in (\omega, \pi/2)$, there exists $C_\mu > 0$ such that the \emph{resolvent estimate} 	$\modulus{\zeta} \norm{ (\zeta - T)^{-1}} \leq C_{\mu}$ holds for all  $\zeta \in \C\setminus \OSec{\mu}$.
\end{enumerate}
\begin{figure}[ht]
\begin{centering}
\includegraphics{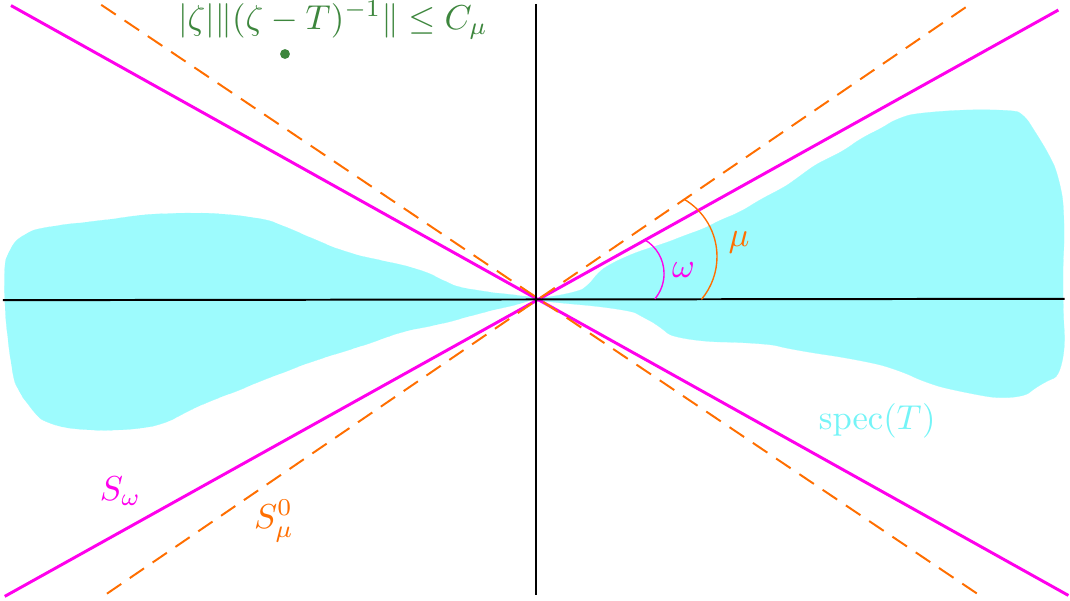}
\end{centering}
\end{figure}
Operators that are $\omega$-bisectorial for which the half plane with non-positive real part $\set{\zeta \in \C: \Re \zeta < 0}$ is free of spectrum are called $\omega$-sectorial.
Every self-adjoint operator is $0$-bisectorial and non-negative self-adjoint operators are $0$-sectorial.

We remark that sectorial operators arise naturally in the study of holomorphic semigroup theory. 
In fact, in \cite{Kato}, there is a notion of sectoriality which further includes the condition on the numerical range, namely that $\Re \inprod{Tu,u}  \geq 0$.
In modern language, these are called \emph{Kato sectorial operators}.

Bisectorial operators automatically admit a rudimentary functional calculus.
Namely, fixing $\mu \in (\omega, \pi/2)$, by $\Psi(\OSec{\mu})$, let us denote the set of holomorphic $\psi: \OSec{\mu} \to \C$ satisfying the estimate: there exists $\alpha > 0$ and $C > 0$ such that
$$\modulus{\psi(\zeta)} \leq C  \max\set { \modulus{\zeta}^\alpha, \modulus{\zeta}^{-\alpha}}$$
for all $\zeta \in \OSec{\mu}$. 
Note that this is reminiscent of the condition \eqref{Eq:Bpassdec}, and in light of this, functions in $\Psi(\OSec{\mu})$ can be thought of as complex valued functions that can be used to define operator adapted band-pass filters. 
For instance,  
$$\zeta \mapsto \frac{\zeta}{1 + \zeta^2}$$
is an example of a $\Psi(\OSec{\mu})$ function. 

Now, a rudimentary functional calculus for $T$ can be defined for each $\psi \in \Psi(\OSec{\mu})$ via the expression
\begin{equation}
\label{Eq:RD}
\psi(T)u := \frac{1}{2 \pi \imath} \oint_{\gamma} \psi(\zeta)(\zeta - T)^{-1}u\ d\zeta,
\end{equation}
where $\gamma$ is the unbounded contour which runs through zero and infinity as
\begin{multline*}
$$\gamma = \set{ t e^{\imath \nu}:\ \infty > t > 0} 
	+ \set{ -te^{-\imath \nu}:\ 0 < t < \infty} \\
	+ \set{ te^{\imath \nu}:\ \infty > t > 0}
	+ \set{ te^{-\imath \nu}:\ 0 < t < \infty}$$
\end{multline*}
with $\nu \in (\omega, \mu)$.
The reason we say it is rudimentary is because via this definition, it is not possible to construct $1(T)$.
However, for $\psi \in \Psi(\OSec{\mu})$, the integral in \eqref{Eq:RD} converges absolutely in the strong operator topology and defines a bounded operator satisfying: there exists a constant $C_{\psi,T}$ (possibly dependent on $\psi$ and $T$) such that 
$$ \norm{\psi(T)u} \leq C_{\psi,T} \norm{u}.$$
Moreover, the bisectoriality of $T$ means that the kernel of $T$ is is accessible from the imaginary axis which, excluding the origin, is contained in the resolvent set with resolvent estimates \ref{Def:Sec:RE}. This induces a decomposition of $\Hil$ as:  
\begin{equation} 
\label{Eq:Split}
\Hil = \nul(T) \oplus \close{\ran(T)}.
\end{equation}
Note that in general, this is a topological sum and not an orthogonal one.
As aforementioned, self-adjoint $T$ are $0$-bisectorial and in this case, this sum is indeed orthogonal.

Despite the ability to create such a functional calculus with many salient features, in applications, it is often necessary to control the norm $C_{\psi,T}$. 
If there exists a universal constant $C > 0$ such that $C_{\psi,T} \leq C\norm{\psi}_{\infty}$, i.e.,
\begin{equation}
\label{Eq:Hinfty}
\norm{\psi(T)u} \leq C \norm{\psi}_{\infty} \norm{u} 
\end{equation}
for all $\psi \in \Psi(\OSec{\mu})$, then we say that $T$ has an \emph{$\Hinfty$ functional calculus}.
This nomenclature will soon be justified.
Let us define a larger class of functions $\Hinfty(\OSec{\mu})$, where we say that $f \in \Hinfty(\OSec{\mu})$ if: 
\begin{enumerate}[(i)] 
\item $f: \OSec{\mu} \union \set{0} \to \C$ is bounded,
\item $f\rest{\OSec{\mu}}: \OSec{\mu} \to \C$ is holomorphic, and  
\item $\norm{f}_\infty < \infty$.
\end{enumerate} 
Note that $f$ is allowed to jump across zero.
In particular, $\sgn(\zeta)  := \chi^+(\zeta) - \chi^-(\zeta)$, where 
\begin{equation} 
\label{Eq:Chipm} 
\chi^{\pm}(\zeta) := \begin{cases} 1 & \pm \Re \zeta > 0, \\  
			0 &\text{otherwise},
	\end{cases}
\end{equation} 
is a  $\Hinfty(\OSec{\mu})$ class function.

Before we extend the rudimentary functional calculus to $\Hinfty(\OSec{\mu})$ functions, let us note the following lemma due to McIntosh. 
\begin{lemma}[McIntosh convergence lemma]
\label{Lem:McIntosh}
If $T$ has an $\Hinfty$ functional calculus, i.e. satisfies \eqref{Eq:Hinfty}, then for each $f\in \Hinfty(\OSec{\mu})$, there exists $f_n \in \Psi(\OSec{\mu})$ such that $f_n \to f$ uniformly on compact subsets of $\OSec{\mu}$ and for each $u \in \Hil$, the sequence $\set{f_n(T)u}$ is Cauchy. 
\end{lemma}
This lemma first appeared in \cite{Mc86} (page 218) by McIntosh.
Armed with this, and writing $\proj_{\nul(T), \close{\ran(T)}}:\Hil \to \nul(T)$ for the projector with kernel $\close{\ran(T)}$ via the splitting \eqref{Eq:Split}, we define
\begin{equation*}
f(T) u := f(0) \proj_{\nul(T), \close{\ran(T)}} u + \lim_{n\to \infty} f_n(T)u.
\end{equation*}
This defines a bounded operator $f(T)$ with norm
$$ \norm{f(T)} \leq C \norm{f}_\infty.$$
It is readily checked that this limit is independent of the particular sequence $\set{f_n(T)u}$ as guaranteed by Lemma \ref{Lem:McIntosh}.
Moreover, this justifies the language ``$\Hinfty$ functional calculus''.
In this setting, by definition, we are using the space $\Psi(\OSec{\mu})$ as a dense subset of $\Hinfty(\OSec{\mu})$.
In a private communication, Alan McIntosh alluded to the construction of the $\Hinfty$ functional calculus as reminiscent to the construction of the Fourier transform in $\Lp{2}(\R^n)$, where the $\Psi(\OSec{\mu})$ can be thought of as ``complex valued Schwartz class functions''.

Bisectorial operators and their accompanying functional calculi are very naturally motivated.
A quintessential example is the first-order factorisation of the \emph{Kato square root problem}.
Let $x \mapsto A(x) \in \Lp{\infty}(\mathrm{SymMat}(\R^n))$, real self-adjoint at almost-every $x$ and $a: \R\to [0,\infty]$ a measurable function.
\nomenclature{$\mathrm{SymMat}(\R^n)$}{Symmetric matrices on $\R^n$.}
We assume that both $A$ and $a$ are uniformly bounded above and below. 
Define the operator:
\begin{equation}
\label{Eq:DivForm} 
\Div_{A,a} u := -a \divv A \nabla,
\end{equation} 
which is seen to be densely-defined and closed via divergence-form methods. 
Of interest is the following nonlinear perturbation estimate:
\begin{equation}
\label{Eq:Qpert}
\norm{\e^{-t L_{A_1,a_1}} - \e^{-t L_{A_2,a_2}}}_{\Lp{2} \to \Lp{2}} \lesssim \norm{A_1 - A_2}_{\Lp{\infty}} + \norm{a_1 - a_2}_{\Lp{\infty}}.
\end{equation}
Note that here, only the  $\Lp{\infty}$ difference of the coefficients appear in the bound - there are no first-order differences.

In order to analyse this operator, the authors of  \cite{AKMc} consider the operator
$$\Pi_{B,b} := \begin{pmatrix} 0 & -b \divv B \\ \nabla & 0 \end{pmatrix},$$
where the matrix-valued map $x \mapsto B(x)$ and the function $x \mapsto b(x)$ are both complex valued and bounded. 
Although a moot point, we remark that $B(x)$ need not be symmetric.
We also assume that there exists $\kappa > 0$ such that $B$ and $b$ satisfy the following accretivity condition:
$$ \Re\inprod{B u, u} \geq \kappa \norm{u}^2\quad\text{and}\quad \Re \inprod{av,v} \geq \kappa \norm{v}^2.$$
The operator $\Pi_{B,b}$ is an $\omega$-bisectorial  operator on $\Hil = \Lp{2}(\R^n)\oplus^\perp \Lp{2}(\R^{2n})$ (since $\nabla: \SobH{1}(\R^n) \subset \Lp{2}(\R^n) \to \Lp{2}(\R^{2n})$) with $\omega < \pi/2$ dependent only on $\kappa$ and the $\Lp{\infty}$ bound on $B$.  
\nomenclature{$\SobH{k}(\cE)$}{Sobolev space of $k$-th order in $\Lp{2}$ on a space $\cE$.}
Then, note that its square is precisely 
$$\Pi_{B,b}^2 = \begin{pmatrix} -b\divv B \nabla  & 0 \\ 0 & -\nabla b \divv B \end{pmatrix},$$
and if $b = a$ and $B = A$, the operator \eqref{Eq:DivForm} is accessed by projecting to the first entry.

Now, let us assume that $\Pi_{B,b}$ has an $\Hinfty$ functional calculus. 
Then, as aforementioned, $f(\Pi_{B,b})$ defined for $f \in \Hinfty(\OSec{\mu})$.
Moreover, we have that $(B',b') \mapsto f(\Pi_{B',b'})$ is holomorphic in a small ball around $(B,b)$ in the $\Lp{\infty}$ topology whose radius is determined by the $\kappa$ and the $\Lp{\infty}$ bounds on $B$ and $b$. 
Consequently, we obtain the Lipschitz estimate 
\vspace{-0.3cm}
\begin{equation*}
\norm{ f(\Pi_{B,b}) - f(\Pi_{B',b'}) } \lesssim \norm{B - B'}_{\Lp{\infty}} + \norm{b - b'}_{\Lp{\infty}},
\end{equation*} 
for $(B',b')$ in this ball.
In particular, on setting $f_t(\zeta) = \e^{-t \zeta^2}$  for $t > 0$, we obtain
$$
f_t(\Pi_{B,b}) = \e^{-t\Pi_{B,b}^2} = \begin{pmatrix} \e^{t b\divv B \nabla}  & 0 \\ 0 & \e^{t\nabla b \divv B} \end{pmatrix}.$$
Then, for $u \in \Lp{2}(\R^n)$,
\begin{multline*} 
\norm{(\e^{-t L_{A_1,a_1}} - \e^{-t L_{A_2,a_2}})u} = \snorm{\cbrac{\e^{-t \Pi_{B_1,b_1}^2} - \e^{-t \Pi_{B_2,b_2}^2}}  \begin{pmatrix} u \\ 0 \end{pmatrix}}   \\ 
	\lesssim (\norm{A_1 - A_2}_{\infty} + \norm{a_1 - a_2}_\infty) \norm{u}.
\end{multline*}
This is precisely the desired estimate in \eqref{Eq:Qpert}.

Although it is possible to obtain such deep consequences in a conceptual manner via the $\Hinfty$ functional calculus, proving that an operator actually possesses this calculus is not often an easy task. 
The criterion \eqref{Eq:Hinfty} is hardly quantitative.
For that, let us return back to the general situation of the $\omega$-bisectorial operator $T$ on a Hilbert space $\Hil$.
At this level of generality, a remarkable theorem due to McIntosh comes to our aid.
It harnesses the perspective that we took in \S\ref{Sec:Discrete} of the reconstruction of the signal in norm as a criterion for the  $\Hinfty$ functional calculus.
The theorem we allude to is the following, although we note that in its original form, it was proved for $\omega$-sectorial operators.
\begin{theorem}[McIntosh's Theorem]
\label{Thm:McIntosh}
Let $T$ be an  $\omega$-bisectorial operator, $\omega \in [0, \pi/2)$,  on a Hilbert space $\Hil$.
Then,  $\Hil = \nul(T) \oplus \close{\ran(T)}$ and for $\mu \in (\omega, \pi/2)$, the following are equivalent:
\begin{enumerate}[(i)]
\item there exists $C > 0$ such that for all $\psi \in \Psi(\OSec{\mu})$, 
	$$\norm{\psi(T)}_{\Hil \to \Hil} \leq C \norm{\psi}_\infty,$$
\item there exists $\psi \in \Psi(\OSec{\mu})$, not identically zero on either sector, such that 
	\begin{equation}
	\label{Qest} 
	\norm{u}_{\psi,T}^2 := \int_{0}^\infty \norm{\psi(tT)u}^2\ \frac{dt}{t} \simeq \norm{u}^2
	\end{equation}
	for all $u \in \close{\ran(T)}$,
\item for all $\psi \in \Psi(\OSec{\mu})$ not identically zero on either sector, \eqref{Qest} holds
	for all $u \in \close{\ran(T)}$.
\end{enumerate}
\end{theorem}
See \cite{Mc86} (page 221) for a proof of this theorem. 
It is also a central theme of the book \cite{Haase} by Haase. 

Of significance are the estimates \eqref{Qest}, which are aptly called \emph{quadratic estimates} in the literature.
This is an incredibly useful quantitative criterion that provides a method with which to assert that $T$ has an $\Hinfty$ functional calculus.
In the same vein as before,  $\psi(tT)$ can be through of as a band-pass filter, localised about frequencies with real part $\frac{1}{t}$.
In fact, the terminology ``band-pass filter'' was initially introduced to the author by Andreas Rosén in a conversation regarding the conceptual foundations of quadratic estimates for bisectorial operators.
Our earlier exposition for operators with discrete spectrum in \S\ref{Sec:Discrete} was borne out of the desire to verify this perspective in a more concrete setting.
Moreover, note that at this level of generality, there is no restriction on type of spectrum that the operator may admit.
In particular, it may have both residue and continuous spectrum in addition to point spectrum.

In the classical setting, via the Fourier transform, the celebrated Plancherel's theorem asserts that the square norm  integral of a  function is equal to the square norm integral of its Fourier transform. 
Via the band-pass filter perspective, it can be seen that \eqref{Qest} is also of this form.
Similarly, in the classical setting, the Calderón reproducing formulas allow the reconstruction of a function through the use of two radial convolutions.
Such reproducing formulas also exist at this level of generality through the $\Hinfty$  functional calculus using specific pairs of $\Psi(\OSec{\mu})$ functions.
This is more or less immediate from the quadratic estimate perspective of the $\Hinfty$ functional calculus.

It is also worth pointing out  a Hilbert space geometric interpretation of quadratic estimates. 
The estimates \eqref{Qest} actually define a new, but comparable Hilbert space norm on $\close{\ran(T)}$. 
Therefore, the ability to assert an $\Hinfty$ functional calculus is equivalent to saying that there are many comparable norms on $\Hil$ which are adapted to the operator.
Given recent developments in synthetic notions of curvature in metric spaces, this interpretation might lend a way of accessing the $\Hinfty$ functional calculus via metric space geometric methods in future.
For a more detailed description of recent developments in metric geometry, see the survey article \cite{Lott} by Lott and references therein. 

For a more detailed treatment of functional calculi for bisectorial and sectorial operators,  we cite the survey paper \cite{ADMc} by Albrecht, Duong, and McIntosh, as well as \cite{Haase}. 
There are also Banach space versions of the $\Hinfty$ functional calculus, initially developed by Cowling, Doust, McIntosh and Yagi in \cite{CDMcY}. 
As this is a grossly inexhaustive list, we urge the reader to consult the references therein.
\section{Geometry} 

\subsection{Kato square root problem in the smooth setting}
\label{Sec:GKato} 
We begin our presentation of a selection of geometric problems resolved through the use of functional calculus by introducing the Kato square root problem in the smooth manifold setting.
To that end, let  $(\cM,\mg)$ be a smooth, complete Riemannian manifold. 
We emphasise that it  need not be compact.
The Laplacian with respect to $\mg$ is then given by
\begin{equation*} 
\Lap_\mg := -\divv_{\mg} \nabla,
\end{equation*}
where $\divv_{\mg} := -\nabla^{\ast,\mg}$ and $\nabla: \SobH{1}(\cM) \subset \Lp{2}(\cM) \to \Lp{2}(\cotanb \cM)$. 
Let $\mh$ be another metric with $C \geq 1$ satisfying 
\begin{equation}
\label{Eq:Gcomp} 
C^{-1} \modulus{u}_{\mg(x)} \leq \modulus{u}_{\mh(x)} \leq C  \modulus{u}_{\mg(x)}
\end{equation}
for all $u \in \tanb_x \cM$.
\nomenclature{$\tanb_x \cM$}{Tangent space of $\cM$ at $x$.} 
\nomenclature{$\cotanb_x \cM$}{Cotangent space of $\cM$ at $x$.} 
Let us define an extended distance metric on the set of all such Riemannian metrics via 
$$ \met_{\cM}(\mg,\mh) = \begin{cases} \infty 	& \text{no $C \geq 1$ exists satisfying \eqref{Eq:Gcomp}}, \\
					\inf\set{\log{C}: C\ \text{satisfies}\ \eqref{Eq:Gcomp}} & \text{otherwise}.
			\end{cases}$$
\nomenclature{$\met_{\cM}(\mg,\mh)$}{Extended distance metric measuring the bounded distance between two metric tensors $\mg$ and $\mh$ on $\cM$.}%
It is easily verified that there exists 
$A \in \Ck{\infty} \intersect \Lp{\infty}(\mathrm{Sym} (\cotanb \cM \tensor \tanb \cM))$ such that
for all $u, v \in \tanb_x \cM$, 
\nomenclature{$\tanb \cM$}{Tangent bundle of $\cM$.}
\nomenclature{$\cotanb \cM$}{Cotangent bundle of $\cM$.}
\begin{equation*}
\mh_{x}(u,v) = \mg_{x}(A(x)u,v)\quad\text{and}\quad
d\mu_{\mh}(x) = \theta(x)\ d\mu_{\mg}(x),
\end{equation*}
where $\theta(x) := \sqrt{\det A(x)}$.
Then,
\begin{equation}
\label{Eq:PLap}
\Lap_{\mh} = - \theta^{-1} \divv_{\mg}A \theta \nabla.
\end{equation}
See \cite{BRough} where these claims are demonstrated in detail.

As before, let us write
$$ \Pi_{B,b,\mg} = \begin{pmatrix} 0 & -\theta^{-1}b \divv_{\mg}B\theta \\ \nabla & 0 \end{pmatrix},$$
for accretive $B \in \Lp{\infty}(\cotanb \cM \tensor \tanb \cM)$ and accretive  $b \in \Lp{\infty}(\cM)$.
Then, as before, from operator theoretic considerations alone, it is possible to prove that $\Pi_{B,b,\mg}$ is an $\omega$-bisectorial operator. 
If $\Pi_{B,b,\mg}$ has an $\Hinfty$ functional calculus, via a similar calculation as before, whenever $\mh$ is another metric with with $\met_{\cM}(g,h) \leq C$, where $C$ is fixed and depends on $(\cM,\mg)$, we obtain that 
\begin{equation}
\label{Eq:Gheat}
\norm{\e^{-t \Lap_{\mg}} - \e^{-t \Lap_{\mh}}}_{\Lp{2}\to\Lp{2}} \lesssim \met_{\cM}(\mg,\mh).
\end{equation}
In fact, more generally, for any $f \in \Hinfty(S_{\mu}^o)$ with $\mu \in (\omega, \pi/2)$, we obtain that 
\begin{equation}
\label{Eq:Gfun}
\norm{f(\Lap_{\mg}) - f(\Lap_{\mh})}_{\Lp{2}\to \Lp{2}} \lesssim \met_{\cM}(\mg,\mh).
\end{equation}
The significance here is that in \eqref{Eq:Gheat}, it is effectively the comparison of two heat kernels with respect to two different metrics, and this estimate says that they are continuous in an $\Lp{\infty}$ sense under the change of metric.
A direct, operator theoretic approach would yield that it is close in the $\Lp{\infty}$ topology in first-derivatives of the change of metric.

The estimate \eqref{Eq:Gfun} is obtained from the generalisation of the Kato square root problem to a geometric setting.
They were initially established in the compact setting in \cite{AKMc} in 2002.
More precisely, the authors of \cite{AKMc}  prove the following. 
\begin{theorem}[Theorem 7.1 in \cite{AKMc}]
Let $(\cM,\mg)$ be a smooth compact Riemannian manifold and suppose that $B \in \Lp{\infty}(\tanb \cM \tensor \cotanb \cM)$ and $b \in \Lp{\infty}(\cM)$ such that
there exists a $\kappa > 0$ satisfying 
$$ \Re \mg_x(B(x)v,v) \geq \kappa \modulus{v}_{\mg(x)}\quad\text{and}\quad \Re b(x) \geq \kappa$$
for almost-every $x \in \cM$ and for all $v \in \cotanb_x \cM$.
Then, the operator $\Pi_{B,b,\mg}$ has an $\Hinfty$ functional calculus with the bound in the estimate depending only on the $\Lp{\infty}$ norms of $B$, $b$ as well as $(\cM,\mg)$ and $\kappa$. 
\end{theorem}
 
A similar problem, but with $\nabla$ replaced by 
$$S = \begin{pmatrix} \nabla \\ \iden \end{pmatrix}: \Lp{2}(\cM) \to \Lp{2}(\cotanb\cM) \oplus \Lp{2}(\cM),$$
with domain $\dom(S) = \SobH{1}(\cM)$, was considered by Morris in 2012 in \cite{Morris}. 
For him, $\cM \subset \R^{n+k}$ is an embedded submanifold.
With this, he considers  the inhomogeneous operator
$$\tilde{\Pi}_{\tilde{B},\tilde{b},\mg} = \begin{pmatrix} 0 & \tilde{b} \adj{S} \tilde{B} \\
					S & 0
			\end{pmatrix},$$
where here $\tilde{B} \in \Lp{\infty}(\cM) \oplus \Lp{\infty}(\tanb \cM \tensor \cotanb\cM)$ and $\mg$ is the induced metric on $\cM$ from $\R^{n+k}$.
Note that the nomenclature ``inhomogeneous'' arises from the fact that the first entry in the square of this matrix realises an operator
$$ \Div_{\tilde{B},\tilde{b}} = \tilde{b}\cbrac{- \divv \tilde{B}_{11} \nabla - \divv \tilde{B}_{10} + \tilde{B}_{01} \nabla + \tilde{B}_{00}},$$
from which its easy to see that this is a perturbation of the Laplacian through the introduction of lower order terms.
He proves the following.
\begin{theorem}[Theorem 1.1 in \cite{Morris}]
Let $\cM \subset \R^{n+k}$ be a smooth, complete Riemannian submanifold and suppose that there exists $C < \infty$ such that the second fundamental form satisfies the bound $\modulus{\SecondFF(\cM)(x)} < C$.
Moreover, suppose there exists $\kappa_1, \kappa_2 > 0$ such that 
$$\Re \inprod{\tilde{B}S u, S u} \geq \kappa_1 \norm{u}^2_{\SobH{1}}\quad\text{and}\quad\Re \inprod{\tilde{b}v,v} \geq \kappa_2 \norm{v}^2.$$
Then, $\tilde{\Pi}_{\tilde{B},\tilde{b},\mg}$ has an $\Hinfty$ functional calculus with the implicit constant dependent on the $\Lp{\infty}$ norms of $\tilde{B}$ and $\tilde{b}$ as well as $C$, $\kappa_1$, $\kappa_2$ and $\cM$.
\end{theorem}

This theorem was subsequently improved to the intrinsic setting in \cite{BMc} by McIntosh and the author to obtain the following.
\begin{theorem}[Theorem 1.1 in \cite{BMc}] 
\label{Thm:InKato}
Let $(\cM,\mg)$ be a smooth, complete Riemannian manifold and suppose there exists $C < \infty$ and $\kappa > 0$ with  
$$\modulus{\Ric_\mg} \leq C\quad\text{and}\quad  \inj(\cM,\mg) \geq \kappa > 0.$$
Suppose that the following ellipticity condition holds: there exist $\kappa_1, \kappa_2 > 0$ such that 
$$\Re \inprod{\tilde{B}S u, S u} \geq \kappa_1 \norm{u}^2_{\SobH{1}}\quad\text{and}\quad\Re \inprod{\tilde{b}v,v} \geq \kappa_2 \norm{v}^2.$$
for all $u \in \Lp{2}(\cM)$ and $v \in \SobH{1}(\cM)$.
Then $\tilde{\Pi}_{\tilde{B},\tilde{b},\mg}$ has an $\Hinfty$ functional calculus with the implicit constant dependent on the $\Lp{\infty}$ norms of $\tilde{B}$ and $\tilde{b}$ as well as $C$, $\kappa$, $\kappa_1$, $\kappa_2$ and $(\cM,\mg)$.
\end{theorem}

A remarkable geometric consequence of these theorems is that, given any $\mgt$ with $\met_{\cM}(\mg, \mgt) < \infty$, we are able to obtain
the statement corresponding to  \eqref{Eq:Gfun} for $(\cM,\mgt)$ in a sufficiently small $\Lp{\infty}$ neighbourhood of $\mgt$.
The size of this neighbourhood, in addition to the geometric assumptions on $(\cM,\mg)$, now also depends on $\met_{\cM}(\mg,\mgt)$.
This observation was exploited in \cite{BRough} to study the necessity of geometric assumptions needed to prove these theorems. 
In particular, it was shown that the curvature and injectivity radius bound assumptions in  Theorem \ref{Thm:InKato} can be dropped.
This indicates that there is ample room for improving the methods used in the proofs.

As a concluding remark, we note that the last theorem, Theorem \ref{Thm:InKato},  is a special case of a more general Kato square root theorem on vector bundles in \cite{BMc}. 
We have abstained from presenting the most general result, which requires some technicalities to facilitate its description, in order to make this exposition more accessible.

\subsection{Rough metrics and a geometric flow tangential to the Ricci flow}
\label{Sec:GFlow}
In \cite{BRough}, a notion of a \emph{rough metric} was defined.
We first note that in order to facilitate the definition of this object, it is important to note that the notion of ``measurable'' and ``null measure'' on a manifold, even when it is noncompact, is independent of a particular choice of Riemannian metric. 
That is to say, we can say that a set is measurable if the intersection of the set with a coordinate chart is measurable against the pullback measure in that chart. 
It is easily verified that for any smooth or even continuous metric tensor $\mg$, the notion of $\mu_\mg$-measurable agrees with this notion of measurable.
Similarly, a null measure set is where the intersected set is null via the pullback measure in each chart.
This allows us to talk about measurable sections of tensorfields, and we define a rough metric to be a symmetric, measurable $(2,0)$-tensorfield for which at each $x \in \cM$, there exists a coordinate chart $(\psi_x, U_x)$ with a constant $C_x = C(U_x) \geq 1$ satisfying
$$ C_x^{-1} \modulus{u}_{\mg(y)} \leq \modulus{u}_{\psi_x^\ast \delta(y)} \leq C_x \modulus{u}_{\mg(y)}$$
for almost-every $y \in U_x$.

It is unclear how to obtain a distance for this metric tensor since the length function is not well-defined.
However, it induces a Borel-regular measure $\mu_\mg$, defined via the usual expression $d\mu_{\mg}(x) = \sqrt{\det \mg(x)}\ d\Leb(x)$. 
As in the smooth setting, this is readily checked to be coordinate independent. 
As before, the notions of measurable  that we remarked were independent of a particular Riemannian structure, despite the lack of regularity of $\mg$, agrees with $\mu_\mg$.
Moreover, a set is of null measure in our sense if and only if  it is null with respect to the $\mu_\mg$ measure.

In this context,  it is again possible to ask the question \eqref{Eq:Gheat}, and more generally \eqref{Eq:Gfun}, for $(\cM,\mg)$.
As before, this follows from establishing an $\Hinfty$ functional  calculus for $\Pi_{B,b,\mg}$.
From the functional calculus perspective, it can also be shown that if $\Pi_{B,b,\mg}$ has an $\Hinfty$ functional calculus, then so does $\Pi_{B,b,\mh}$ for any rough metric $\mh$ that satisfies \eqref{Eq:Gcomp}.
This shows that such metric tensors are geometric invariances of the Kato square root problem, which further leads to questions about the existence of geometries which are  counterexamples to this problem, and consequently, $\Hinfty$ functional calculus.

There are many natural geometries that are rough metrics.
Take, for instance, a smooth metric $\mh$ and a lipeomorphism (locally bi-Lipschitz map) $\phi: \cM \to \cM$, and consider the metric $\mg = \phi^\ast \mh$.
The metric $\mg$ is a rough metric, and in general it will only have measurable coefficients since the definition of this metric involves derivatives of the Lipschitz map $\phi$.
Note that even in $\R^2$, given a null measure set, there exists a Lipschitz function for which the non-differentiable points are contained in this null set. 
This demonstrates the existence of crude and non-trivial rough metrics even in the limited context of lipeomorphic pullbacks alone. 
In particular, the set of singularities can be a dense subset.
From a geometric point of view, rough metrics are useful when a manifold can be given a smooth differentiable structure so that singularities can be represented purely in terms of a lack of regularity of the metric tensor. 
Such is the case, for instance, for a conical singularity with spherical cross sections.

Beyond being an object of interest for the reasons we have just outlined, rough metrics have become useful in the study of a geometric flow with non-smooth initial data.
In \cite{GM}, Gigli and Mantegazza defined a geometric flow on a smooth compact manifold $(\cM,\mg)$, which they ultimately generalise to measure-metric spaces satisfying a certain synthetic Ricci curvature bound.
The flow in the smooth setting is obtained by solving for $\phi_{t,x,v}$ in the \emph{continuity equation}
\begin{equation}
\label{Def:GMC}
-\divv_{\mg,y} \hk^\mg_t(x,y) \conn \phi_{t,x,v}(y) = \extd_{x}(\hk^\mg_t(x,y))(v)
\end{equation}
for each fixed $x \in \cM$, $t > 0$.
Here, $\hk^\mg_t$ is the heat kernel of $\Lap_\mg$, the operator $\divv_{\mg,y}$ denotes the divergence operator acting on the variable $y$, the vector $v \in \tanb_x \cM$, and $\extd_x(\hk^\mg_t(x,y))(v)$ is the directional derivative of $\hk^\mg_t(x,y)$ in the variable $x$ in the direction $v$.
They define a new family of metrics
evolving in time by the expression
\begin{equation}
\label{Def:GM}
\mg_t(x)(u,v) = \int_{\cM} \mg(y)(\conn \phi_{t,x,u}(y), \conn \phi_{t,x,v}(y))\ \hk^\mg_t(x,y)\ d\mu_\mg(y).
\end{equation}
Moreover, they demonstrate that this flow satisfies
$$\partial_t \mg_t(\dot{\gamma}(s),\dot{\gamma}(s))\rest{t = 0} = -2 \Ric_{\mg}(\dot{\gamma}(s),\dot{\gamma}(s))$$ 
for almost-every $s$ along $\mg$-geodesics $\gamma$.
That is, this flow $t \mapsto \mg_t$ is \emph{tangential} (in this weak sense) to the Ricci flow.
While the authors of \cite{GM} do not touch upon questions of regularity in their paper, they show that $x \mapsto \mg_t(x)$ remains smooth for all $t > 0$.
It is of interest is to understand smoothing or non-smoothing properties of this flow in measure-metric settings.
However, such a task is difficult in the general measure-metric space world  given that we lack even a basic yardstick to make sense of reasonable notions of  regularity. 
However, between the smooth Riemannian manifold setting and the metric space world, there is an entire jungle of non-smooth metric tensors against smooth differentiable structures. 
This is a fertile playground in which we can attempt to address regularity questions.
 
Given that \eqref{Def:GMC} is in divergence form, it can certainly be understood in a weak sense.
Starting with an initial rough metric, the best we can expect for the regularity of the heat kernel is $\Ck{\alpha}$-regularity via the methods of Nash-Moser-de Giorgi.
This is precisely what is done in \cite{BCont}. 
However, to understand the continuity of $x \mapsto \mg_t(x)$ on an open region where the heat kernel improves to $\Ck{1}$, we are forced with a non-smooth perturbation problem with measurable coefficients.
Note from \eqref{Eq:PLap},  by choosing any smooth metric, equation \eqref{Def:GM} is equivalent to solving a divergence form equation with bounded measurable coefficients coming from the regularity of the initial metric $\mg$, regardless of whether $\hk^\mg_t$ improves in regularity in an open region. 
In fact, in \cite{BCont}, via this method, the regularity question is reduced to showing the Kato square root problem for rough metrics. 
There, it is shown that the Kato square root problem can be solved for \emph{any} rough metric on a compact manifold.
More precisely, the following is proved in \cite{BCont}.
\begin{theorem}[Theorem 3.4 and Theorem 4.3 in \cite{BCont}] 
\label{Thm:Kato}
On a compact manifold $\cM$ with a rough metric $\mg$,  the operator $\Pi_{B,b,\mg}$ admits an $\Hinfty$ functional calculus.
Consequently, if  $\emptyset \neq \cN \subset \cM$ is an open set where the initial heat kernel satisfies $\hk^\mg_t \in \Ck{1}(\cN^2)$,
$\mg_t$ as defined by \eqref{Def:GM} exists on $\cN$ and it is continuous.
\end{theorem}

As far as the author is aware, this is the first  and only known instance where the Kato square root problem has been applied to obtain regularity properties of a geometric flow.

\subsection{Elliptic boundary value problems}
Let $(\cM,\mg)$ be a smooth, compact  Riemannian manifold with smooth boundary $\Sigma := \partial \cM$.
Moreover, let $(\cE,\mh^\cE), (\cF, \mh^\cF) \to \cM$ be Hermitian vector bundles and $\Dir:\Ck{\infty}(\cM; \cE) \to \Ck{\infty}(\cM; \cF)$ a first-order elliptic differential operator.
\nomenclature{$\Ck{k}(\cM; \cE)$}{Smooth sections $\cM \to \cE$.}
Here, elliptic means that the principal symbol $\sym_{\Dir}(x,\xi):\cE_x \to \cF_x$ is invertible.
\nomenclature{$\sym_{\Dir}(x,\xi)$}{Principal symbol of $\Dir$ in the co-direction $\xi$ at $x$.}
By $\Ck[c]{k}(\cM; \cE)$, let us denote the $u \in \Ck{k}(\cM; \cE)$ such that $\spt u$ is compact.
\nomenclature{$\Ck[c]{k}(\cM; \cE)$}{Compactly supported smooth sections $\cM \to \cE$ up to the boundary (if it exists).}
In particular, we allow for the possibility that $\spt u \intersect \Sigma \neq \emptyset$. 
On the other hand, by $\Ck[cc]{k}(\cM; \cE)$, we mean $u \in \Ck[c]{k}(\cM; \cE)$ such that $\spt u \subset \interior \cM$, i.e. $\spt u \intersect \Sigma = \emptyset$.
\nomenclature{$\Ck[cc]{k}(\cM; \cE)$}{Compactly supported smooth sections $\cM \to \cE$  on the interior.}
\nomenclature{$\interior \cM$}{Interior of $\cM$.}  
The restrictions $\cE\rest{\Sigma}$ and $\cF\rest{\Sigma}$ define smooth bundles on $\Sigma$.
We remark that our discussion in this section also carries over to the noncompact setting (but still with compact boundary), but there, further assumptions need to be placed on the operator $\Dir$.
We avoid conducting a conversation at this level of generality as it adds an unnecessary layer of complexity that detracts from the key features of this example. 

By local considerations alone, there exists a formal adjoint $\Dir^\dagger: \Ck{\infty}(\cM; \cF) \to \Ck{\infty}(\cM; \cE)$. 
On writing $\Dir_{\cc}^\dagger = \Dir^\dagger$ with $\dom(\Dir_{\cc}^\dagger) = \Ck[cc]{\infty}(\cM; \cF)$, we obtain a maximal operator $\Dir_{\max}$ and a minimal operator $\Dir_{\min}$ by
$$ \Dir_{\max} := \adj{(\Dir_{\cc}^\dagger)}\quad\text{and}\quad  \Dir_{\min} := \close{\Dir_{\cc}}.$$
It is easy to see that both are closed and densely-defined operators.

In order study boundary value problems for first-order differential operators, it is important to be able to define the \emph{boundary restriction map} $u \mapsto u\rest{\Sigma}$, initially defined on $\Ck[c]{\infty}(\cM; \cE)$, on the space $\dom(\Dir_{\max})$. 
Moreover, it is essential that we understand the topology of the range of this map in order to obtain a bounded surjection.
By abstract considerations, noting $\dom(\Dir_{\min})$ is the kernel of this extended boundary restriction map, it is possible to show that this topology is unique.
This allows for a complete description of all closed extensions of $\dom(D_{\min})$ as closed subspaces of the range of the boundary restriction map.
An important motivation comes from index theory. 
As Atiyah-Patodi-Singer demonstrated in \cite{APS,APS1,APS2,APS3}, it is essential to be able to deal with non-local boundary conditions to phrase index theorems.
It is also desirable to be able to understand how the index behaves under perturbations of boundary conditions. 
The panoramic perspective afforded through understanding the boundary restriction map on $\dom(D_{\max})$ is necessary to consider such questions.
 
In \cite{BB}, Bär and Ballmann demonstrate how this can by done by computing the topology of the range of the boundary trace map on $\dom(D_{\max})$  via an associated elliptic first-order differential operator $\DirA$ on $\Lp{2}(\Sigma; \cE)$. 
This adapted operator on the boundary has principal symbol
\begin{equation} 
\label{Eq:BSymb}
 \sym_{\DirA}(x,\xi) = \sym_{\Dir}(x, \tau(x))^{-1} \comp \sym_{\Dir}(x,\xi),
\end{equation} 
where $\xi \in \cotanb_x \Sigma$ and   $x\mapsto \tau(x)$ is an interior pointing co-vectorfield on $\Sigma$.
There are many operators that satisfy this condition, but in their treatise, Bär and Ballmann assume that  the symbol $\sym_{\DirA}(x,\xi)$ is skew-adjoint. 
In this case, there are self-adjoint operators, unique up to the addition of a zeroth order operator, that satisfy \eqref{Eq:BSymb}.
For such an operator, the range 
$$\set{u\rest{\Sigma}: u \in \dom(\Dir_{\max})} =  \checkH(\DirA) := \ind{(-\infty,0)}(A)\SobH{\frac{1}{2}}(\Sigma; \cE) \oplus \ind{[0,\infty)}(\DirA)\SobH{-\frac{1}{2}}(\Sigma; \cE),$$
and $u \mapsto u\rest{\Sigma}$ is a bounded surjection in this topology.
Boundary conditions are precisely closed subspaces $B \subset \checkH(A)$ and the associated closed operator $\Dir_{B,\max}$ to $B$ is $\Dir$ with domain 
$$ \dom(\Dir_{B,\max}) = \set{ u \in \dom(\Dir_{\max}): u\rest{\Sigma} \in B}.$$
Since the operator $\DirA$ is self-adjoint, elliptic, first-order, and the boundary $\Sigma$ is compact, it only admits discrete spectrum with orthogonal eigenspaces. Therefore, the analysis in \cite{BB} is carried out in the spirit of the Fourier series, heavily exploiting the orthogonality of eigenspaces.
Although this condition on $\DirA$ seems restrictive, their results are very general and can be applied liberally to a wide class of operators, including all Dirac type operators.

A quintessential example is the operator $\SDir$, the Atiyah-Singer Dirac operator. 
This induces an operator $\DirA$ such that $\sym_{\DirA}(x,\tau(x))$ anti-commutes with $\DirA$.
Then, the boundary condition 
$$B_{\text{APS}} = \chi_{(-\infty,0)}(\DirA)\SobH{\frac{1}{2}}(\Sigma; \cE)$$
is precisely the famed Atiyah-Patodi-Singer boundary condition used in the proof of their index theorem.
In this situation, the operator $\SDir_{B_{\text{APS}}}$ is self-adjoint.
Given $r \in \R$, we can also consider generalised APS conditions:
$$B_{\text{gAPS}}(r) = \chi_{(-\infty,r)}(A)\SobH{\frac{1}{2}}(\Sigma; \cE).$$
When $r \neq 0$, the operator $\SDir_{B_{\text{APS}}(r)}$ is no longer self-adjoint. 

A more exotic object of interest is the \emph{Rarita-Schwinger} operator. 
Again, the manifold $\cM$ is assumed to be  Spin and  this operator is defined as follows. 
Let $\cW := \cotanb\cM \tensor \Spinors\cM$ where $\Spinors \cM$ is the Spin bundle and define the map 
$\iota: \Spinors\cM \to \cW$ via  
$$\iota(\psi) := -\frac{1}{n} \sum_{j=1}^{n} e_j^\ast \tensor (e_j \cdot \psi).$$ 
Then, let
$$ \Spinors_{\frac{3}{2}}\cM  := \iota(\Spinors\cM)^\perp$$
so that 
$\cW = \iota(\Spinors\cM) \stackrel{\perp}{\oplus} \Spinors_{\frac{3}{2}}\cM.$
There, we obtain an  induced twisted Dirac operator $\SDir_{\cW}$ which can be written as the operator matrix
$$ \SDir_{\cW} = \begin{pmatrix} \iota \SDir & \adj{\mathrm{T}} \\ \mathrm{T} & \SDirRS \end{pmatrix}.$$
The operator $\SDirRS$ is the Rarita-Schwinger operator. 
The principal symbol of any adapted boundary operator $\DirA_{\text{RS}}$ to $\SDirRS$, given by  $\sym_{\DirA_{\text{RS}}}(x,\xi) = \sym_{\SDirRS}(x,\tau(x))^{-1} \comp \sym_{\SDirRS}(x,\xi)$,  fails to be skew-symmetric in general dimensions.
This is one operator that falls  outside of the scope of the framework in \cite{BB}.

In  \cite{BBan}, Bär and the author demonstrate that any differential operator $\DirA$ obtained from a first-order elliptic differential operator $\Dir$ satisfying \eqref{Eq:BSymb} is, up to the addition of a real constant, an invertible $\omega$-bisectorial operator.
In this case,  via the work of Grubb in \cite{Grubb}, they obtain the spectral projectors $\chi^{\pm}(\DirA)$ to the left and right of the complex plane as classical pseudo-differential operators of order zero.
With the aid of this, the operator 
$$\modulus{\DirA} := \DirA\sgn(\DirA)$$
is defined, where $\sgn(\zeta) = \chi^+(\zeta) - \chi^-(\zeta)$ and  where $\chi^{\pm}$ were defined in \eqref{Eq:Chipm}.
By the boundedness of $\chi^{\pm}(\DirA)$ and commutativity with $\DirA$,  it is clear that $\dom(\modulus{\DirA}) = \dom(\DirA)$. 
Moreover, the operator $\modulus{\DirA}$ is $\omega$-sectorial and invertible.
As before, define 
$$\checkH(\DirA) := \chi^-(\DirA)\SobH{\frac{1}{2}}(\Sigma; \cE) \oplus \chi^+(\DirA)\SobH{-\frac{1}{2}}(\Sigma; \cE).$$
Unfortunately, at this level of generality, the analysis is far more subtle than \cite{BB}.
For instance, even the rudimentary task of demonstrating the boundedness of $u \mapsto u\rest{\Sigma}: \dom(\Dir_{\max}) \to \checkH(\DirA)$ can no longer be conducted as in \cite{BB} since the generalised eigenspaces associated to the spectrum of $\DirA$ may no longer be orthogonal. 
It turns out, however, that the analysis precisely reduces to establishing that $\modulus{\DirA}$ has an $\Hinfty$ functional calculus.

The criterion that McIntosh's Theorem \ref{Thm:McIntosh} provides is useful as it allows one to move between the $\Hinfty$ functional calculus and quadratic estimates.
In the analysis in \cite{BBan}, it is really the quadratic estimates that are needed. 
These estimates are accessed via establishing the $\Hinfty$ functional calculus by other means. 
Namely, by a more recent work due to Auscher, Nahmod and McIntosh in \cite{AMcN}, it is possible to assert the $\Hinfty$ functional calculus for an $\omega$-sectorial operator via interpolation methods.
Via their Corollary 5.5 in \cite{AMcN}, it suffices to show that there are two numbers $s, t > 0$ such that  $\dom(\modulus{\DirA}^s) \subset \dom(\modulus{\DirA^\ast}^s)$ and $\dom(\modulus{\DirA^\ast}^t) \subset \dom(\modulus{\DirA}^t)$.
In our case, we can choose $s = t = 1$ from the fact that $\DirA$ is a first-order elliptic differential operator on a smooth bundle $\cE\rest{\Sigma}$ over a compact manifold $\Sigma$.
More precisely, by  pseudo-differential elliptic regularity theory, we obtain  
$$ \dom(A) = \dom(\modulus{A}) = \dom(A^\ast) = \dom(\modulus{A^\ast}) = \SobH{1}(\Sigma; \cE).$$
In fact, Corollary 5.5 in \cite{AMcN} can be thought of as a conduit between pseudo-differential methods and the $\Hinfty$ functional calculus. 
The following is an important result obtained in \cite{BBan}.

\begin{theorem}
\label{Thm:EllMain} 
Let $(\cM,\mg)$ be a compact smooth manifold with compact boundary, and let $\Dir: \Ck{\infty}(\cM; \cE) \to \Ck{\infty}(\cM; \cF)$ be a first-order elliptic differential operator between Hermitian bundles $(\cE,\mh^\cE) \to \cM$ and $(\cF,\mh^{\cF}) \to \cM$. 
Then, the following hold: 
\begin{enumerate}[(i)]
\item $\Ck[c]{\infty}(\cM; \cE)$ is dense in $\dom(\Dir_{\max})$ with respect to the corresponding graph norm.
\item 
\label{Eq:BdyBnd} 
The trace map $\Ck[c]{\infty}(\cM; \cE) \to \Ck{\infty}(\Sigma; \cE)$	given by $u \mapsto u\rest{\Sigma}$ extends uniquely to a surjective bounded linear map
	$\dom( \Dir_{\max}) \to \checkH(\DirA)$.
\item The spaces
\begin{align*} 
	\dom(\Dir_{\max}) &\intersect \SobH[loc]{k}(\Sigma; \cE)\\ 
		&= \set{ u \in \dom(\Dir_{\max}): Du \in \SobH[loc]{k-1}(\cM; \cE)\ \text{and}\ u\rest{\Sigma} \in \SobH{\frac{k}{2}}(\Sigma; \cE)}.
\label{Eq:MaxDom}
\end{align*}
\item The $\Lp{2}$ inner product $\inprod{\cdot,\cdot}$ extends to a perfect paring $\inprod{\cdot,\cdot}_{\checkH(A) \times \checkH(-A^\ast)}$ between $\checkH(A)$ and $\checkH(-A^\ast)$.
For all $u \in \dom(\Dir_{\max})$ and $v \in \dom((\Dir^\dagger)_{\max})$,
\begin{equation}
\inprod{\Dir_{\max} u, v}_{\Lp{2}(\cM; \cF)} - \inprod{u, (\Dir^\dagger)_{\max}v}_{\Lp{2}(\cM; \cE)} = -\inprod{u\rest{\Sigma}, \sym_{\Dir}(\cdot,\tau)^\ast v\rest{\Sigma}}_{\checkH(A) \times \checkH(-A^\ast)}.
\end{equation}
\end{enumerate}
The corresponding statements hold for $\Dir^\dagger$ with the roles $\cE$ and $\cF$ interchanged and on replacing $A$  by $\tilde{A}$, a boundary adapted operator to $\Dir^\dagger$.
\end{theorem}

To illustrate through an example how the $\Hinfty$ functional calculus appears in these proofs, let us consider the following simple calculation.
As we have aforementioned, these methods can be made to work for noncompact manifolds but with compact boundary.
For the purpose of the calculation which we are about to perform, it is beneficial for us to consider $\cM = \R_+ \times \Sigma$  for some closed manifold $\Sigma$. 
\nomenclature{$\R_+$}{The set $[0, \infty)$.}
Despite that it is noncompact, it affords us with a greater simplicity than choosing a compact example.
Moreover, further assume that $\Dir = \partial_t + \DirA$ over some bundle $\cE = \cF$ with $\DirA$ invertible $\omega$-bisectorial first-order differential operator.
In this setting, the operator $\modulus{\DirA}$ has an $\Hinfty$ functional calculus.

Let $\eta \in \Ck[c]{\infty}(\R)$ be a cutoff function such that $\eta = 1$ on $[0, 1]$ and $0$ on $[2,\infty)$.
Then, for $v \in \Ck[c]{\infty}(\Sigma)$, define 
\begin{equation}
\label{Eq:Ext} 
(\sE v)(x,t) =  \eta(t) (\e^{-t\modulus{\DirA}}v)(x).
\end{equation}
By the properties of the semigroup, we have that $(\sE v)\rest{\Sigma}(x) = (\sE v)(x,0) = v(x)$.
Our goal is to show that 
\begin{equation} 
\label{Eq:Req}
\norm{\sE v}_{\Dir} \simeq \norm{\Dir \sE v}_{\Lp{2}(\cM;\cE)} + \norm{\sE v}_{\Lp{2}(\cM;\cE)}  \lesssim \norm{v}_{\checkH(\DirA)},
\end{equation}
which is a crucial step in proving \ref{Eq:BdyBnd} of Theorem \ref{Thm:EllMain}.
So, first note that 
\begin{equation}
\label{Eq:HEComp} 
\begin{aligned} 
&\partial_t (\sE v)(t,\cdot) =  \eta'(t) \e^{-t\modulus{\DirA}} v + \eta(t) \modulus{\DirA} \e^{-t\modulus{\DirA}}v \\
&\DirA (\sE v)(t,\cdot) = \eta(t) \DirA \e^{-t\modulus{\DirA}}v.
\end{aligned}
\end{equation}
Moreover, write $v_{\pm} = \chi^{\pm}(\DirA)v$ so that $v = v_- + v_+$.
This enables us to consider the cases $v_{\pm}$ separately from linearity of \eqref{Eq:Ext} and the boundedness of the projectors $\chi^{\pm}(\DirA)$. 

Since $\modulus{\DirA}u = \DirA u$, from \eqref{Eq:HEComp}, it follows that 
$$ \Dir (\sE v_+)(t, \cdot) = \eta'(t) \e^{-t\modulus{\DirA}}v_+.$$
Thus, it follows easily from the fact that $\norm{\eta'}_{\Lp{\infty}} < \infty$ that
\begin{align*}
\norm{\Dir \sE v_+}^2_{\Lp{2}(\cM;\cE)} &= \int_{0}^\infty \norm{\eta'(t) \e^{-t\modulus{\DirA}}v_+}^2_{\Lp{2}(\Sigma;\cE)}\ dt  \\
	&\lesssim \int_{0}^\infty \norm{t^{\frac{1}{2}}\modulus{\DirA}^{\frac{1}{2}} \e^{-t\modulus{\DirA}} \modulus{\DirA}^{-\frac{1}{2}} v_+}^2_{\Lp{2}(\Sigma;\cE)} \ \frac{dt}{t}.
\end{align*}
On setting $\psi(\zeta) = \zeta^{\frac{1}{2}}\e^{-\zeta}$, we have that $\psi \in \Psi(S_{\mu}^o)$ for $\mu \in (0, \pi)$ and  
$$
\psi(t\modulus{\DirA}) = t^{\frac{1}{2}}\modulus{\DirA}^{\frac{1}{2}} \e^{-t\modulus{\DirA}}.$$ 
Therefore, since we assert that $\modulus{\DirA}$ has an $\Hinfty$ functional calculus, via McIntosh's Theorem \ref{Thm:McIntosh}, we obtain that
\begin{equation}
\label{Eq:Plus}
\int_{0}^\infty \norm{t^{\frac{1}{2}}\modulus{\DirA}^{\frac{1}{2}} \e^{-t\modulus{\DirA}} \modulus{\DirA}^{-\frac{1}{2}}v_+}^2_{\Lp{2}(\Sigma;\cE)}\ \frac{dt}{t} 
	\lesssim \norm{\modulus{\DirA}^{-\frac{1}{2}}v_+}^2 \simeq \norm{v_+}_{\SobH{-\frac{1}{2}}(\Sigma;\cE)}^2.
\end{equation}
Note that  the ultimate equivalence in this calculation follows from elliptic regularity theory. 

For the remaining case, we have that $\modulus{\DirA}v_- = -\DirA v_-$ and from \eqref{Eq:HEComp},
$$ \Dir (\sE v_-)(t, \cdot) = (\eta'(t) - 2\eta(t)\modulus{\DirA})  \e^{-t\modulus{\DirA}}v_-.$$
Therefore, 
$$
\norm{\Dir (\sE v_-)}_{\Lp{2}(\cM;\cE)} \lesssim \int_{0}^\infty \eta'(t)^2 \norm{e^{-t\modulus{\DirA}}v_-}^2_{\Lp{2}(\Sigma;\cE)} \ dt +  \int_0^\infty 4\eta(t)^2 \norm{\modulus{\DirA}e^{-t\modulus{\DirA}}v_-}^2_{\Lp{2}(\Sigma;\cE)} \ dt.$$
The first term is bounded similar to \eqref{Eq:Plus}, but for the second term,
\begin{equation} 
\label{Eq:Minus}
\begin{aligned}
\int_0^\infty \eta(t)^2 \norm{\modulus{\DirA}\e^{-t\modulus{\DirA}}v_-}^2\ dt
	&\lesssim \int_0^\infty \norm{t^{\frac{1}{2}}\modulus{\DirA}^{\frac{1}{2}} \e^{-t \modulus{\DirA}} \modulus{\DirA}^{\frac{1}{2}}v_-}^2 \ \frac{dt}{t}\\
	&\lesssim \norm{\modulus{\DirA}^{\frac{1}{2}}v_-}^2 \simeq \norm{v_-}_{\SobH{\frac{1}{2}}(\Sigma; \cE)}^2.
\end{aligned}
\end{equation}
The penultimate and crucial inequality again relies upon the $\Hinfty$ functional calculus of $\modulus{\DirA}$. 
On combining \eqref{Eq:Plus} and \eqref{Eq:Minus} yields \eqref{Eq:Req}.

As a concluding remark, let us  emphasise that  despite the discreteness of the spectrum of $\DirA$, in this more general setting, this does not allow us to reduce the problem to considering individual generalised eigenspaces of $\DirA$.
One issue is their potential non-orthogonality. 
This apparent inconvenience is, in fact, a boon, as it has inspired the development of methods that are robust enough to be applied to even more general settings.
For instance, the functional calculus perspective developed in understanding the general elliptic operator situation for compact boundary can, in fact, be applied to adapted operators that are symmetric on the boundary when the boundary is non-compact.
Under a set of mild hypotheses, this yields a self-adjoint operator. 
The spectrum may now contain residue and continuous blobs in addition to points, but the $\Hinfty$ functional calculus methods are completely insensitive and are able to accommodate for this situation.
This calculus exists for these operators since every self-adjoint operator enjoys an $\Hinfty$ functional calculus.
Through careful analysis, it is possible to prove a theorem, similar in spirit to Theorem \ref{Thm:EllMain}, in the non-compact boundary case.

\subsection{Spectral flows and Riesz continuity for the Atiyah-Singer Dirac operator} 
\label{Sec:SFlow}
In \S\ref{Sec:GKato}, we demonstrated that an $\Hinfty$ functional calculus for the operator $\Pi_{B,b,\mg}$ yields Lipschitz estimates \eqref{Eq:Gfun} for small perturbations of the coefficients $(B,b)$ in the $\Lp{\infty}$ topology.
This is obtained via establishing a much stronger result, namely that $(A,a) \mapsto f(\Pi_{A,a,\mg})$ is holomorphic in an
$\Lp{\infty}$ neighbourhood of $(B,b)$.
The proof of this exploits the fact that $\Pi_{B,b,\mg}$ admits the structure
\begin{equation}
\label{Eq:KatoStr} 
\Pi_{B,b,\mg} = \Gamma + \tilde{B}_1 \adj{\Gamma} \tilde{B}_2,
\end{equation}
where $\Gamma$ is a nilpotent, closed, densely-defined operator, and the $\tilde{B}_1$ and $\tilde{B_2}$ appropriately encode the coefficients $B$ and $b$. 
In particular, this affords us with a Hodge-decomposition in $\Lp{2}$. 
Many desirable consequences follow from \eqref{Eq:KatoStr} and it is unclear how these could be obtained without such structure.
 
A geometric operator that is also of physical significance is the Atiyah-Singer Dirac operator $\SDir_{\mg}$ over the Spin bundle $\Spinors \cM$ when the manifold $\cM$ is Spin with a Riemannian structure given by $\mg$. 
There are some natural questions that one can ask.
The first is in the case that  $\cM$ has compact boundary $\Sigma = \partial \cM$, and where $\cB$ is a self-adjoint boundary condition for $\SDir_{\mg}$, yielding an operator $\SDir_{\cB,\mg}$. It is interesting to know what happens to the operator when we perturb the boundary condition $\cB$ in an appropriate sense.
A primary complication  is that the domain $\dom(\SDir_{\cB,\mg})$ moves when we perturb $\cB$.
In fact, in a private communication with Alan Carey, he revealed that the late Krzysztof Wojciechowski was present at the Australian National University in 2004, and he was excited by the first-order factorisation of the Kato square root problem  as well as the corresponding problem for the Hodge Dirac operator in \cite{AKMc} which naturally accounts for moving domains. 
He wondered whether these methods maybe applicable to study perturbations of boundary conditions for Dirac operators.
His ambition was to understand the perturbation of such operators in the Riesz topology, which has connections to index theory and spectral flows.

A second related question, perhaps simpler in nature, is on a manifold without boundary. 
Here, it is desirable to to understand what happens to $\SDir_{\mg}$ under perturbation of the metric $\mg$. 
An added complication in this scenario that the Spin bundle $\Spinors\cM$, outside of a set of exceptional perturbations, itself might change when the metric moves.

Wojciechowski's excitement is indeed justified.
It is possible to consider these problems via methods motivated by the resolution of the Kato square root problem.
However, unlike this problem which can be seen via \eqref{Eq:KatoStr}, it is unclear the kind of structural features that $\SDir$ might enjoy.
Consequently, we are forced to consider the Lipschitz estimate directly, and for this, it is actually necessary to go beyond the regime of a single bisectorial operator and instead consider the functional calculus of two self-adjoint operators.
We will provide a brief exposé of the details later.
First, let us present some results from \cite{BMcR}  by Rosén, McIntosh and the author where the metric perturbation question was addressed, and \cite{BR}  by Rosén and the author where the boundary condition perturbation problem was studied. 

Let us consider the question of metric perturbation.
As we have aforementioned, the bundle itself can change here. 
However, if we are interested in understanding spectral properties, we are able to pull back through a similarity transform using a unitary operator.
This allows us to compare the spectra of two operators despite the fact that they may act on two different Spin bundles.
More precisely, if we have another metric $\mh$ on $\cM$, we can construct a unitary map $U:(\tanb\cM, \mg) \to (\tanb\cM, \mh)$, which naturally extends to a map between Spin bundles. 
In what is to follow, let 
$$\OSec{\omega,\sigma}:=\set{x+iy: y^2<\tan^2\omega x^2 + \sigma^2}.$$
The metric perturbation theorem is then the following. 

\begin{theorem}[Theorem 3.1 in \cite{BMcR}]
\label{Thm:SFMet}
Let $\cM$ be a smooth Spin manifold (without boundary) with a smooth, complete metric $\mg$ and Levi-Civita  connection $\conn^\mg$. 
Fixing  constant $C_{\mh} > 0$, let  $\mh$ be a $\Ck{0,1}$ metric with $\met_{\cM}(\mg, \mh) \leq 1$ satisfying:
\begin{enumerate}[(i)]
\item 
\label{Hyp:MainAppFirst}
\label{Hyp:Inj}
there exists $\kappa > 0$ such that $\inj(\cM,\mg) \geq \kappa$,
\item 
\label{Hyp:Curv}
there exists $C_{R} > 0$ such that 
	$\modulus{\Ric_\mg} \leq C_R$ and $\modulus{\conn^\mg \Ric_\mg} \leq C_R$, 
\item 
\label{Hyp:MainAppLast} 
	$\modulus{\conn^\mg \mh} \leq C_{\mh}$
	almost-everywhere.
\end{enumerate}

Then, for $\omega \in (0, \pi/2)$, $\sigma > 0$, whenever $f \in \Hinfty(\OSec{\omega,\sigma})$, 
we have the perturbation estimate
$$
\norm{f(\SDir_\mg) - f(\spin{\U}^{-1}\SDir_\mh \spin{\U})}_{\Lp{2} \to \Lp{2}}	
	\lesssim \norm{f}_{\infty} \met_M(\mg,\mh).$$
The implicit constant depends on $\dim\cM$ and the constants appearing in \ref{Hyp:MainAppFirst}-\ref{Hyp:MainAppLast}. 
\end{theorem}

Returning to the original question posed by Wojciechowski, when the manifold $\cM$ has boundary $\Sigma$, we can consider perturbations of local boundary conditions for the Atiyah-Singer Dirac operator using the framework  in \cite{BB}. 
A \emph{local boundary condition} is a space  
$$\cB = \SobH{\frac{1}{2}}(\cE)\quad \text{with}\quad \cE \subset \Spinors\Sigma = \Spinors \cM \rest{\Sigma},$$ 
where $\cE$ is a smooth subbundle.   
The operator $\SDir$ with boundary condition $\cB$, denoted
$\SDir_{\cB}$, is the operator with domain 
$$\dom(\SDir_{\cB}) = \set{ \phi \in \dom(\SDir_{\max}): u\rest{\Sigma} \in \cB}.$$

To describe our perturbation result, we require two additional conditions on the local boundary condition $\cB$: 
\begin{enumerate}[(i)]
\item Self-adjointness, which by \S3.5 in \cite{BB} occurs
if and only if $\sym_{\SDir}(\on^\flat)$ maps the $\Lp{2}$
closure of $\cB$ onto its orthogonal complement.
Here $\on^\flat$ the unit normal interior covectorfield.

\item $\SDir$-ellipticity, which is defined in terms
of a self-adjoint boundary operator $\BDir$ adapted
to $\SDir$ with principal symbol 
$\sym_{\BDir}(\xi) = \sym_{\SDir}(\on^\flat)^{-1} \comp \sym_{\SDir}(\xi)$, 
and for which the operator 
$$\pi_{\cB} - \ind{[0, \infty)}(\BDir):\Lp{2}(\Spinors\Sigma) \to \Lp{2}(\Spinors\Sigma)$$ 
is a Fredholm operator. Here, 
$\pi_{\cB}: \Lp{2}(\Spinors\Sigma) \to \cB$ 
is projection induced from the fibrewise orthogonal
projection $\pi_{\cE}: \Spinors\Sigma \to \cE$, and
$\ind{[0,\infty)}(\BDir)$ is the projection
onto the positive spectrum of the operator $\BDir$ (see Theorem 3.15 in \cite{BB}).
\end{enumerate} 

For two local boundary conditions 
$\cB$ and $\tilde{\cB}$,   following \S2 in Chapter IV in \cite{Kato},
we define the \emph{$\Lp{\infty}$-gap}  between the subspaces $\cB$ and $\tilde{\cB}$ as
$$ \gap_\infty(\cB, \tilde{\cB}) := \norm{\gap(\cE_x, \tilde{\cE}_x)}_{\Lp{\infty}(\Sigma)} := \sup_{x\in \Sigma} \modulus{\pi_{\cE}(x) - \pi_{\tilde{\cE}}(x)},$$
where $\pi_{\cE}$ and $\pi_{\tilde{\cE}}$ are the
orthogonal projections from $\Spinors\Sigma$ to $\cE$ and 
$\tilde{\cE}$ respectively.
We let  
$$\norm{\cB}_{\Lips} := \sup_{x\in \Sigma} \modulus{\conn \pi_{\cE}(x)},$$   
and similarly for $\tilde{\cB}$. 
For a set $Z \subset \cM$ which is a neighbourhood of $\Sigma$, given  $r > 0$, we write $Z_r = \set{x \in \cM: \met_\mg(x, Z) < r}$. 
By $Z_r \disunion Z_r$, we denote the double of $Z_r$ by pasting along $\Sigma$.

\begin{theorem}[Theorem 3.1 in \cite{BR}]
\label{Thm:SFBdy}
Let $(\cM,\mg)$ be a smooth, Spin manifold with smooth, compact boundary $\Sigma = \bnd\cM$
 that is complete as a metric space  and suppose:
\begin{enumerate}[(i)]
\item  there exists a precompact open neighbourhood $Z$ of $\Sigma$  and $\kappa > 0$ such that $\inj(\cM\setminus Z,\mg) > \kappa$, 
\item there exists $C_R < \infty$ such that $\modulus{ \Ric_\mg} \leq C_R$ 
	and $\modulus{\conn \Ric_\mg}  \leq C_R$ on $\cM \setminus Z$, and
\item  fix any smooth metric $\mg_Z$ on the double $Z_4 \disunion Z_4$ obtained by pasting along $\Sigma$ and $C_Z < \infty$ and let $\kappa_Z > 0$ 
	with $\modulus{\Ric_{\mg_Z}} \leq C_Z$ and $\inj(Z_2 \disunion Z_2, \mg_Z) \geq \kappa_Z$. 
\end{enumerate}  

Fixing $C_B < \infty$, let  $\cB$ and $\tilde{\cB}$ be two local self-adjoint $\SDir$-elliptic boundary conditions which satisfy: 
\begin{enumerate}[(i)]
\item[(iv)] $\norm{\cB}_{\Lips} + \norm{\tilde{\cB}}_{\Lips} \leq C_B$, and
\item[(v)] $\SDir$-ellipticity constants of orders $1$ and $2$ for $\cB$ in a given compact neighbourhood $K$ of the boundary.
\end{enumerate}  
Then, for $\omega \in (0, \pi/2)$ and $\sigma > 0$, 
whenever we have $f \in \Hinfty(\OSec{\omega, \sigma})$, we have the perturbation
estimate
$$ \norm{f(\SDir_\cB) - f(\SDir_{\tilde{\cB}})}_{\Lp{2}\to \Lp{2}} \lesssim \norm{f}_\infty \gap_\infty(\tilde{\cB}, \cB).$$
The implicit constant depends on $\dim\cM$ and the constants appearing in (i)-(v).
\end{theorem}

The sole motivation for accounting for the function class $\Hinfty(\OSec{\omega,\sigma})$ is that it includes 
$$\zeta \mapsto \frac{\zeta}{\sqrt{1 + \zeta^2}}.$$
This allows us to measure the distance between the two operators in the Riesz topology. 
This is a desirable metric to work with in attempting to understand the spectral flow for these operators as it connects better with topological properties of these operators, as seen through K-theory. 
This was observed in \cite{AS69} by Atiyah and Singer for bounded self-adjoint Fredholm operators.
Moreover, note that for $\cM$ compact, the geometric hypotheses in both Theorem \ref{Thm:SFMet} and Theorem \ref{Thm:SFBdy} are automatically satisfied.
A virtue of these results, even when applied to the compact setting,  is that they quantify how these quantities enter the perturbation estimate.

In both these results, it is important to note that the continuity is in an $\Lp{\infty}$ sense.
That is, like for the perturbation estimate in the Kato square root problem, the continuity depends only in an $\Lp{\infty}$ sense of the perturbation, although here, we are also required to perturb in a way in which we are uniformly $\Lp{\infty}$ bounded in the gradient of the metric or the projector defining the local boundary condition.
This is contrast to earlier results where an additive term of $\norm{\nabla^\mg \mh}_{\Lp{\infty}}$ or $\norm{\cB}_{\Lips} + \norm{\tilde{\cB}}$ would appear in the right side of the estimate.

A concrete example of a family of metrics $\mg_{\epsilon}$ on $\R^2$ that fits the hypothesis of Theorem \ref{Thm:SFMet} are: 
$$ \mg_\epsilon(x) = \begin{pmatrix} 1 + \epsilon \sin\cbrac{\frac{\modulus{x}}{\epsilon}} & 0 \\ 
				0 & 1 
		\end{pmatrix}.$$
For $\epsilon > 0$, note that these metrics are not  conformal to the standard Euclidean metric, but they are $\Ck{0,1}$ with the property that  $\norm{\nabla^{\R^2} \mg_{\epsilon}}_{\Lp{\infty}} \leq 1.$

As we have remarked earlier, in contrast to the resolution of the Kato square root problem, the analysis here is not to establish  Lipschitz estimates through holomorphic dependency results for a single operator that encodes the perturbation. 
Instead, we are forced to consider the difference of functional calculus of two different operators.
This is a complicated matter despite the fact that the operators  are self-adjoint.
Denoting the perturbation of either question by operators $\Dir_1$ and $\Dir_2$, and on drawing inspiration from  methods arising from the resolution of the Kato square root problem, we are able to reduce the estimate of  $\norm{ f(\Dir_1) - f(\Dir_2)}_{\Lp{2} \to \Lp{2}}$
to a \emph{local quadratic estimate} of the form 
$$ 
\int_{0}^1 \snorm{ \frac{t\Dir_1}{1 + t^2 \Dir_1^2} X \frac{t\Dir_2}{1 + t^2 \Dir_2^2} u}_{\Lp{2}}\ \frac{dt}{t} \lesssim \norm{X}_{\Lp{\infty}} \norm{u}^2_{\Lp{2}}.$$
The coefficients $X$ either encode the difference between the metrics or the gap metric for the case of boundary condition perturbations.
It is a bounded operator, a multiplication operator, or a special kind of singular integral operator.
This reduction is  obtained on drawing inspiration from McIntosh's Theorem \ref{Thm:McIntosh}, where quadratic estimates are shown to imply (in fact, are equivalent to) the $\Hinfty$ functional calculus.

Note here that the upper limit of the integral is $1$ and not  $\infty$ as in the case of the quadratic estimates \eqref{Qest} in Theorem \ref{Thm:McIntosh}.
This is effectively because functions  in $\Hinfty(\OSec{\omega,\sigma})$ are holomorphic on a neighbourhood of $0$, and on recalling that $\psi(t\Dir)$ can be thought of as a localisation of a signal to the spectrum of $\Dir$  with real part $\frac{1}{t}$, we can expect that there are only contributions from  high frequencies.
\section{Harmonic Analysis}

In our exposition so far, we have already used the phrase ``band-pass filter'' in \S\ref{Sec:FunCalc}. 
This is suggestive, at least superficially, of potential links between functional calculus and harmonic analysis. 
Indeed, a particular significance of the quadratic estimates perspective of the $\Hinfty$ functional calculus is that it is a bridge to the realm of real-variable harmonic analysis from which these estimates can be calculated. 
In the two examples described in  \S\ref{Sec:GKato} and \S\ref{Sec:GFlow}, we have seen that the the goal has been to obtain an $\Hinfty$ functional calculus via quadratic estimates. 
In \S\ref{Sec:SFlow}, despite the fact that this problem is outside of the scope of the $\Hinfty$ functional calculus, there is a reduction to a ``local'' quadratic estimate.
The aim of this section is to illustrate, using the philosophy which emerges from the resolution of the Kato square root problem, of how such an estimate might be computed.

Lets fix a scenario before we flesh out the intended outline. 
Fix $(\cM,\mg)$ to be a complete Riemannian manifold and suppose that $\Dir_1$ and $\Dir_2$ are first-order differential operators on a Hermitian bundle $(\cE,\mh^E) \to \cM$.
For simplicity, assume and that they are self-adjoint in $\Lp{2}(\cM;\cE)$.
We remark that, since we allow for vector bundles, the first-order nature of these operators are not a limitation.
This can be seen in the setting of \S\ref{Sec:GKato}, where the problem is actually second-order for functions.
Here, we have factorised this problem into a first-order system involving functions and vectorfields. 
Even in the situation of $(\cM,\mg) = (\Sph^2,\mg_{\text{round}})$, the two sphere with the round metric, the tangent and the cotangent bundles are not globally trivialisable. This is seen easily from the hairy ball theorem which states that any continuous vectorfield on $\Sph^2$ must vanish at some point.
This simple example illustrates the reason which the vector bundle setting is unavoidable. 

At a first glance, we have seemingly complicated the situation by allowing for vector bundles.
Worse still, a factorisation of an elliptic higher-order problem might result in a non-elliptic first-order system.
However, this tall cost is easily justified.
A feature of a first-order operator $\Dir$ is that the principal symbol, or equivalently the commutator of $\Dir$ with functions, is a multiplication operator. 
Often, such operators are geometric, by which we mean that it arises from the geometry in some natural way. 
Quantitatively, this means that there is some $C > 0$ such that 
\begin{equation} 
\label{Eq:CommBd} 
\modulus{[\Dir, f\iden](x)} \leq C \modulus{\nabla f(x)}.
\end{equation}
This, along with an assumption $f \dom(\Dir) \subset \dom(\Dir)$ for $f \in \Ck[c]{\infty}(\cM)$ and the self-adjointness of $\Dir$, is enough to show that this operator admits \emph{exponential off-diagonal estimates}. 
That is, on fixing a $\psi \in \Psi(\OSec{\mu})$, there exists $C_{\psi} > 0$ and for each $M > 0$, there exists a constant $C_{\Delta,\psi, M} > 0$ so that 
\begin{equation}
\label{Def:OD}
\begin{aligned}
\norm{\ch{E} \psi(t\Dir) (\ch{F}u)}_{\Lp{2}(\cM; \cE)}
	\leq  C_{\Delta,\psi, M}  &\maxx{\frac{\met(E,F)}{t}}^{-M} \times \\ 
		&\qquad\exp\cbrac{-C_{\psi} \frac{\met(E,F)}{t}} \norm{\ch{F}u}_{\Lp{2}(\cM; \cE)}  
\end{aligned}
\end{equation}
for every Borel set $E,\ F \subset \cM$ and $u \in \Lp{2}(\cM; \cE)$.
These estimates can be found in Lemma 5.3 in \cite{CMcM} by Carbonaro, McIntosh and Morris for a specific choice of operator.
A closer inspection of their proof reveals that the only necessary ingredients are the global boundedness of the commutator as in \eqref{Eq:CommBd} and the fact that the domain of the operator is preserved under multiplication by compactly supported smooth functions.

The significance of such an  estimate is that it allows us to localise computations.
However, the functional calculus is built out of holomorphic functions and they admit tails when localised.
In order to ensure convergence over a sum of norms of localisations, an estimate of the form \eqref{Def:OD} is required.
More precisely,  given a disjoint cover $\set{B_j}$, it is desirable to localise by decomposing as 
$$ \norm{\psi(t\Dir)S v} = \sum_{j} \norm{\psi(t \Dir) S v}_{\Lp{2}(B_j)} \leq \sum_{j} \sum_{i} \norm{\chi_{B_j} \psi(t \Dir) S \ind{B_i} v},$$
where  $S$ is typically a  multiplication operator.
It is clear from the  term on the right where off-diagonal estimates assist in the computation.
These estimates can be thought of as a replacement for Gaussian estimates, the bread and butter of the scalar valued second-order world. 
However, as we have seen, for a natural class of first-order operators, exponential off-diagonal decay is obtained gratis. 
Moreover, this is also suggestive that our analysis might be carried out in spaces with exponential growth such as negatively curved manifolds.
We shall touch upon this shortly.

Without  loss of generality, let us assume that we are attempting to prove quadratic estimates of the form 
\begin{equation} 
\label{Eq:Qneed} 
\int_{0}^{\tau} \norm{\QQ_t X \PP_t u}^2\ \frac{dt}{t} \lesssim \norm{X}^2_\infty \norm{u}^2,
\end{equation}
where either $\tau = 1$ or $\tau = \infty$.
Here, the operator 
$$\PP_t = \frac{1}{1 + t^2 \Dir_2^2},$$
whereas the operator $\QQ_t = \psi(t\Dir_1)$ for some $\psi \in \Psi(\OSec{\mu})$.
This satisfies
$$ \int_{0}^1 \norm{\QQ_t u}^2\ \frac{dt}{t} \lesssim \norm{u}^2,$$
as well as the exponential off-diagonal estimates  \eqref{Def:OD}.
The quintessential example of such an operator is 
$$\QQ_t = \frac{t\Dir_1}{1 + t^2 \Dir_1^2}.$$

In order to proceed, we need to also understand how the underlying measure-metric geometry of the manifold $(\cM,\mg)$ affects the ability to perform harmonic analysis.
This means we need to exploit some notion of scale invariance coupled to the geometry in a useful manner.
In Euclidean space, this is captured through dyadic cubes.
This is just a decomposition of $\R^n$ into a grid, where at each scale $j \in \In$, the cubes  of the grid are of length $2^j$ and are a nested subdivision in a natural way of the cubes of scale $j+1$ with length $2^{j+1}$.

Such structures exist in more general settings, and for us,  the most significant result is due to Christ in \cite{Christ}.
There, he constructs a dyadic structure that naturally takes into account the measure-metric structure of $(\cM,\mg)$. 
As observed in \cite{MorrisThesis}, his result is easily generalised to spaces of \emph{exponential volume growth}. More precisely, this means that there exists $c_E \geq 1$, $\kappa, c > 0$
such that 
\begin{equation}
\label{Def:Eloc}
0 < \mu_\mg(\Ball(x,tr)) \leq ct^\kappa \e^{c_E tr} \mu_\mg(\Ball(x,r)) < \infty, 
\end{equation}
for every $t \geq 1,\ r > 0$ and $x\in \cM$.
For instance, if there exists $\eta \in \R$ such that the Ricci curvature satisfies  $\Ric_{\mg} \geq \eta \mg$, then $(\cM,\mg)$ satisfies \eqref{Def:Eloc}. 
Under this assumption, the theorem is as follows. 
\begin{theorem}[Existence of a truncated dyadic structure]
\label{Thm:Dya:Christ}
Suppose that $(\cM,\mg)$ satisfies \eqref{Def:Eloc}.
Then, there  exist countably  many index sets $I_k$,
a countable collection of open subsets
$\set{\Q[\alpha]^k \subset \cM: \alpha \in I_k,\ k \in \Na}$,
points $z_\alpha^k \in \Q[\alpha]^k$ (called the \emph{centre} of $\Q[\alpha]^k$), 
and constants $\delta \in (0,1)$, 
$a_0 > 0$, $\eta > 0$ and $C_1, C_2 < \infty$ satisfying:
\begin{enumerate}[(i)]
\item for all $k \in \Na$,
	$\mu(\cM \setminus \union_{\alpha} \Q[\alpha]^k) = 0$,
\item if $l \geq k$, then either $\Q[\beta]^l \subset \Q[\alpha]^k$
	or $\Q[\beta]^l \intersect \Q[\alpha]^k = \emptyset$,
\item for each $(k,\alpha)$ and each $l < k$
	there exists a unique $\beta$ such that
	$\Q[\alpha]^k \subset \Q[\beta]^l$,
\item $\diam \Q[\alpha]^k < C_1 \delta^k$, 
\item $ \Ball(z_\alpha^k, a_0 \delta^k) \subset \Q[\alpha]^k$,  
\item for all $k, \alpha$ and for all $t > 0$, 
	$\mu\set{x \in \Q[\alpha]^k: d(x, \cM\setminus\Q[\alpha]^k) \leq t \delta^k} \leq C_2 t^\eta \mu(\Q[\alpha]^k).$
\end{enumerate} 
\end{theorem}

As in the Euclidean case, this theorem allows us to obtain a \emph{dyadic grid}, suitability adjusted to the distance and volume structure of $(\cM,\mg)$, beneath some scale $\scale > 0$.
Note this is in contrast to the Euclidean setting, where the dyadic grid is arbitrarily large. 
This is possible more generally in non-compact settings if the measure-metric structure of $(\cM,\mg)$ satisfies the stronger condition of \emph{doubling}, rather than the weaker notion of exponential volume growth that we assume.

This dyadic decomposition is essential to the construction of a \emph{dyadic averaging operator},
$$\Av_t:\Lp{2}(\cM; \cE) \to \Lp{2}(\cM; \cE),$$
uniformly bounded for $t \leq \scale$.
Moreover, an  essential feature of this operator is that, in an appropriate sense, it is ``constant'' for $x \in \Q$ for each dyadic cube $\Q$.
This notion requires some mild bounded geometry assumptions on the bundle and the manifold, and we shall remark on this towards the end.
Via this notion of constancy, along with the off-diagonal decay coupled to the  dyadic decomposition, we can construct a certain \emph{principal part} operator 
$$\Pri_t(x): \cE_x \to \cE_x.$$ 
This is a kind of local smoothing operator and we will soon see that this is a key operator in the analysis.
Using these operators, the underlying philosophy  at the heart of the resolution of the Kato square root problem is to break up the integrand in \eqref{Eq:Qneed} in the following way: 
\begin{equation}
\label{Eq:SFEBreak}
\begin{aligned}
\int_0^1 \norm{\QQ_tS\Ppb_t f}^2\ \dtt
	&\lesssim  \int_0^1 \norm{(\QQ_t - \Pri_t\Av_t) X\PP_t f}^2\ \dtt \\
		&\qquad+ \int_0^1 \norm{\Pri_t\Av_tX(\iden - \PP_t) f}^2\ \dtt 
		+ \int_0^1 \norm{\Pri_t \Av_t X f}^2\ \dtt.
\end{aligned}
\end{equation}

The first term, 
\begin{equation} 
\label{Eq:Poin} 
\int_0^1 \norm{(\QQ_t - \Pri_t\Av_t) X\PP_t f}^2\ \dtt,
\end{equation}
aptly called the principal part term,  is estimated by resorting to a Poincaré inequality, appropriately generalised to the bundle setting to match with the notion of constancy that we have aforementioned in passing.
The second term,
\begin{equation}
\label{Eq:Can} 
\int_0^1 \norm{\Pri_t\Av_tX(\iden - \PP_t) f}^2\ \dtt,
\end{equation}
reduces to exploiting certain cancellation properties of the differential operator, and again, this is typically a mild term to estimate. 

The key purpose of this decomposition is really to reduce the whole problem, through these two mild estimates, to a \emph{local Carleson measure} estimate.
While the first two estimates can be treated by ``soft'' methods, that is, more or less exploiting operator theoretic properties of the problem, it is in the local Carleson measure estimate where recent developments in modern real-variable harmonic analysis become of great importance.

To describe this in further detail, let us provide the definition of a \emph{local Carleson measure}. 
A measure  $\nu$ is called a local Carleson measure on $\cM \times (0,t']$ (for some fixed $t' \in (0, \scale]$) if
$$\norm{\nu}_{\Carl} = \sup_{t \in (0,t']} \sup_{\Q \in \DyQ_t} \
	\frac{\nu(\CBox(\Q))}{\mu(\Q)} < \infty.$$ 
Here,  $\CBox(\Q) = \Q \times (0, \len(\Q))$ is the  \emph{Carleson box} over $\Q$, $\len(\Q)$ is the length of the cube $\Q$, and $\DyQ_t$ are the dyadic cubes whose length $\len(\Q) \sim t$.
The norm $\norm{\nu}_{\Carl}$ is the \emph{local Carleson norm} of $\nu$.
For a local Carleson measure $\nu$, Carleson's inequality yields 
$$\iint_{\cM \times (0,t']} \modulus{\Av_t(x)u(x)}^2\ d\nu(x,t) \lesssim \norm{\nu}_{\Carl} \norm{u}^2$$
for all $u \in \Lp{2}(\cM; \cE)$.
A proof of this, in the case of functions, can be found as Theorem 4.2 in \cite{Morris}, but this proof is equally valid in our vector bundle setting.

Returning back to estimating the remaining term
$$\int_0^1 \norm{\Pri_t \Av_t X f}^2\ \dtt,$$
note that by what we have just mentioned, it suffices to prove that 
$$ d\nu(x,t) = \modulus{\Pri_t(x)}^2\ \frac{d\mu_\mg(x)dt}{t}$$ 
is a local Carleson measure.
This typically involves heavy duty harmonic analysis machinery including non-tangential maximal functions, local $T(b)$ theorems, and other ideas arising from real-variable harmonic analysis.
It is unfortunately beyond the scope of this article to give a description of these methods, but  the books \cite{Stein70, Stein93} by Stein and
\cite{Christ90} by Christ give detailed accounts of these ideas.
The survey article \cite{HMc} touches upon more recent developments of these ideas, particularly in connection with the Kato square root problem.

As we have aforementioned, this procedure for obtaining the quadratic estimates require mild bounded geometry assumptions on both the manifold and the bundle.
For instance, in Theorem \ref{Thm:InKato}, Theorem \ref{Thm:SFMet}, and Theorem \ref{Thm:SFBdy}, this is the reason that Ricci curvature and injectivity radius bounds appear as an assumption.
We have already remarked that a bounded geometry assumption is required for the notion of constancy to which we alluded earlier. 
We also require this for a generalised local Poincaré inequality which is bootstrapped from functions.
These are required to estimate \eqref{Eq:Poin} as well \eqref{Eq:Can}.

It is unclear how necessary these assumptions are to carry out quadratic estimates.
These assumptions exist, in part, due to the fact that these methods are a generalisation of real-variable harmonic analysis machinery that has been developed in the Euclidean setting, where geometrically speaking, the bundles are flat and trivial.
However, it is reasonable to expect,  given that we now have an awareness of the geometric issues that arise in the analysis, that these methods may be better adapted  to account for the geometry of the problem. 
Investigations along these lines may allow for the bounded geometry notion to be further weakened, or perhaps dispensed with entirely. 
While this is expected to be a gross undertaking that is technical in nature,  this would be a worthwhile task as it would likely reveal deeper links between geometry and harmonic analysis through the medium of functional calculus.

\printnomenclature[2.5cm]

\bibliographystyle{amsplain}
\def\cprime{$'$}
\providecommand{\bysame}{\leavevmode\hbox to3em{\hrulefill}\thinspace}
\providecommand{\MR}{\relax\ifhmode\unskip\space\fi MR }
\providecommand{\MRhref}[2]{%
  \href{http://www.ams.org/mathscinet-getitem?mr=#1}{#2}
}
\providecommand{\href}[2]{#2}

\setlength{\parskip}{0mm}

\end{document}